\documentclass[11 pt, reqno]{article}

\usepackage{amsmath,amsfonts,amssymb,amsthm}
\usepackage[colorlinks=true,hyperindex=true]{hyperref}
\usepackage{cancel,bbm}
\usepackage{mathabx} 
\usepackage{comment}
\usepackage{enumitem}
\usepackage{cite}
\usepackage{comment}

\usepackage{chngcntr}
\counterwithin*{equation}{section}

\usepackage[left=3cm, right=3cm, top=3 cm, bottom= 3 cm]{geometry}
\usepackage{color}
\usepackage{fancyhdr}
\usepackage{latexsym}

\usepackage{bm}

\newtheorem{ccounter}{ccounter}[section]
\newtheorem{thm}[ccounter]{Theorem}
\newtheorem{lem}[ccounter]{Lemma}
\newtheorem{cor}[ccounter]{Corollary}

\newtheorem{prop}[ccounter]{Proposition}
\newtheorem{ass}[ccounter]{Assumption}
\newtheorem{ex}[ccounter]{Example}

\newtheorem{definition}[ccounter]{Definition}
\newtheorem{rem}[ccounter]{Remark}

\def\bet{\begin{thm}}
\def\eet{\end{thm}}
\def\bel{\begin{lem}}
\def\eel{\end{lem}}
\def\bas{\begin{ass}}
\def\eas{\end{ass}}
\def\bec{\begin{cor}}
\def\eec{\end{cor}}
\def\bep{\begin{prop}}
\def\eep{\end{prop}}
\def\beq{\begin{equation}}
\def\eeq{\end{equation}}
\def\proof{\noindent {\bf Proof.}\ \ }
\def\bea{\begin{equation*}}
\def\eea{\end{equation*}}

\def\bex{\begin{ex}}
\def\eex{\end{ex}}

\def\nn{\mathbb{N}}
\def\rr{\mathbb{R}}

\def\cc{\mathbb{C}}
\def\1{\boldsymbol{1}}
\let\Im\relax
\DeclareMathOperator{\Im}{Im}

\def\e{\mathrm{e}}
\def\i{\mathrm{i}}
\def\del{\partial}
\def\d{\mathrm{d}}
\def\eps{\varepsilon}
\renewcommand\leq\varleq
\renewcommand\geq\vargeq
\def\ee{\mathrm{E}}

\def\O{\mathcal{O}}

\def\ee{\mathbb{E}}

\def\pp{\mathbb{P}}

\def\D{\mathcal{D}}

\def\D{\mathcal{D}}

\def\be{\mathbf{e}}

\def\E{\mathcal{E}}

\newcommand{\corr}{\color{red}}
\newcommand{\cobb}{\color{blue}}
\newcommand{\nc}{\normalcolor}
\def\1{\boldsymbol{1}}
\def\mft{\mathfrak{t}}
\newcommand{\av}{\mathrm{ave}}
\newcommand{\iso}{\mathrm{iso}}
\renewcommand{\le}{\leq}
\renewcommand{\ge}{\geq}
\def\be{{\bf e}}
\def\hata{\hat{a}}
\def\hatb{\hat{b}}

\begin{document}


\begin{table}
\centering

\begin{tabular}{c}

\multicolumn{1}{c}{\parbox{12cm}{\begin{center}\Large{\bf Optimal Delocalization for Non--Hermitian Eigenvectors}\end{center}}}\\
\\
\end{tabular}
\begin{tabular}{ c c c  }
Giorgio Cipolloni
& \phantom{blah} & 
Benjamin Landon
 \\
 & & \\  
 \small{University of Rome Tor Vergata} & & \small{ University of Toronto } \\
 \small{Department of Mathematics} & & \small{Department of Mathematics} \\
 \small{\texttt{cipolloni@axp.mat.uniroma2.it}} & & \small{\texttt{blandon@math.toronto.edu}} \\
  & & \\
\end{tabular}
\\
\begin{tabular}{c}
\multicolumn{1}{c}{\today}\\
\\
\end{tabular}

\begin{tabular}{p{15 cm}}
\small{{\bf Abstract:} We prove an optimal order delocalization estimate for the eigenvectors of general $N \times N$ non-Hermitian matrices $X$:  $\| {\bf v } \|_\infty  \leq C \sqrt{\frac{\log N}{N}}$ with very high probability, for any right or left eigenvector ${\bf v}$ of $X$. This improves upon the previous tightest bound of Rudelson and Vershynin \cite{MR3405592} of $\O ( \frac{ ( \log N)^{9/2}}{ \sqrt{N}})$, and holds under weaker assumptions on the tail of the matrix elements.  In addition to the coordinate basis, our bound holds for the $\ell^\infty$ norm in  any deterministic orthonormal basis.

\phantom{tab} Our result is proven via a dynamical method, by studying the flow of the resolvent of the Hermitization of $X$ and proving local laws on short scales. 

}
\end{tabular}
\end{table}

\section{Introduction}

The delocalization and localization properties of the eigenvectors of large random matrices play an important role in modeling physical phenomena. For $N \times N$ random matrices comprised of Gaussian entries the question of delocalization  has essentially a complete solution: the distribution of the matrix is invariant under rotations and so any $\ell^2$-normalized eigenvector ${\bf v}$ is uniformly distributed on the sphere. Therefore,
\beq \label{eqn:int-v}
\lVert {\bf v} \rVert_\infty \leq C \sqrt{\frac{ \log N}{N}} 
\eeq
with very high probability. On the other hand, one no longer has an explicit distribution for the eigenvectors of more general (e.g., non-Gaussian) random matrices and obtaining an estimate as sharp as \eqref{eqn:int-v} becomes a challenging problem.

Nonetheless, a variety of methods have been developed to tackle delocalization problems for very general classes of random matrices. The theory is quite well-developed for random Hermitian matrices: see for example, \cite{aggarwal2019bulk, erdHos2009local, gotze2020local, knowles2013isotropic, tao2012random, vu2015random} (and also the recent surveys \cite{bourgade2018random, erdHos2017dynamical, o2016eigenvectors}). In fact, in the Hermitian case, even the optimal constant in \eqref{eqn:int-v} is known \cite{benigni2022optimal} for generalized Wigner matrices (mean field random Hermitian matrices with independent entries whose variance matrix is stochastic).

Despite the successes in the Hermitian case, the non-Hermitian theory is underdeveloped. The best result is due to the seminal work of Rudelson and Vershynin  \cite{MR3405592} who developed a geometric approach to delocalization and proved that
\beq
\lVert {\bf v} \rVert_\infty \leq C \frac{ ( \log N)^{9/2}}{\sqrt{N}} ,
\eeq
with very high probability. Their result holds for a general class of mean-field non-Hermitian matrices with sub-Gaussian entries. However, the above estimate falls short of the optimal upper bound \eqref{eqn:int-v}. In the work \cite{nguyen2018normal}, an estimate of optimal order  was obtained for the eigenvalues of small modulus (i.e., eigenvalues $\sigma$ s.t. $| \sigma| \leq C N^{-1/2}$ in our scaling) for matrices $X$ comprised of i.i.d. entries, but this constitutes only a vanishing fraction of the 
 entire spectrum of $X$. 

We will consider  $N \times N$ matrices $X$ with i.i.d. centered entries. Our main result shows that for any $D>0$ there is a $C = C(D)$ so that for any $\ell^2$-normalized left or right eigenvector ${\bf l}$ or ${\bf r }$ of $X$ we have,
\begin{equation}
\label{eq:evdeloc-a1}
\lVert {\bf l}\rVert_\infty+\lVert {\bf r}\rVert_\infty\le C\sqrt{\frac{\log N}{N}},
\end{equation}
with probability at least $1-N^{-D}$. In particular, this is the first estimate of optimal order for general i.i.d. matrices.  Our work improves on all prior results in the following three ways:
\begin{enumerate}[label=(\roman*)]
\item The factor $(\log N)^{9/2}$ has been reduced to the optimal $( \log N)^{1/2}$.
\item We do not require any sub-Gaussian or any sub-exponential tails for the matrix elements, only finiteness of high moments.
\item Our estimates are not limited to the coordinate basis of $\cc^N$, but hold in any deterministic basis.
\end{enumerate}
 In a different direction, delocalization has been proven for more general non-Hermitian random matrices with a variance profile \cite{alt2021spectral}, a deformation \cite{campbell2024spectral}, or even for matrices with certain correlations \cite{alt2021inhomogeneous}, at the price of a worse error $N^\xi/\sqrt{N}$, for any small $\xi>0$. Additionally, the work of Rudelson and Vershynin \cite{MR3405592} holds (under a uniform sub-Gaussian assumption) for matrices of arbitrary variance profile. We expect that our work can easily be extended to models with a large class of variance profiles and deformations  (see the remark at the end of Section~\ref{sec:mainres}).

 We remark here that the upper bound in \eqref{eq:evdeloc-a1} cannot rule out that most of the eigenvector coordinates are $0$. There have been many works proving \emph{lower bounds} for the mass that an eigenvector assigns to sets of indices; this has sometimes been called no gaps delocalization after the work of Rudelson and Vershynin \cite{rudelson2016no} which proved one of the first such results for general matrices. Further work improving on \cite{rudelson2016no} in the case of i.i.d. matrices appears in \cite{lytova2020delocalization,luh2020eigenvector}. 

 There has additionally been significant progress in determining the fluctuations of the individual entries of the eigenvectors of i.i.d. matrices, with the asymptotic Gaussian fluctuations having been obtained in \cite{dubova2024gaussian,osman2025least}. Note, however, that the convergence in distribution alone is, of course, insufficient to obtain the high probability estimates of our work.

\subsection{Methodology}

One of the main technical tools to study $X$ is the following \emph{Hermitization} of $X$, defined by, 
\begin{equation}
H^z :=\left(\begin{matrix}
0 & X-z \\
(X-z)^* & 0
\end{matrix}\right).
\end{equation}
In particular, $z \in \mathrm{spec}(X)$ iff $0 \in \mathrm{spec}(H^z)$ and the corresponding eigenvectors are related in a simple fashion. It therefore suffices to consider the resolvent $G^z(\i\eta):=(H^z-\i\eta)^{-1}$ and rely on the fact that the size of the eigenvectors ${\bf w}_i^z$ of $H^z$ can by estimated by
\[
\big|\langle{{\bf w}_i^z,{\bf x}}\rangle\big|^2\lesssim \eta \cdot \langle {\bf x}, \Im G^z(\i\eta){\bf x}\rangle,
\]
for any $\eta>0$. In the bulk of the spectrum when $|z| \leq 1 - c$, it suffices to choose $\eta = C \frac{\log N}{N}$ and show that $\langle {\bf x}, \Im G^z(\i\eta){\bf x}\rangle\lesssim 1$ with very high probability.\footnote{At the edge of the spectrum the scaling is more delicate and we omit this for the present discussion.} 

The main difficulty is therefore to prove a local law for $\langle {\bf x}, \Im G^z(\i\eta){\bf x}\rangle$ in this very small $\eta$-regime. The standard methodology of proving local laws based on matrix methods typically requires relatively large scales; i.e., that  $\eta\ge N^\xi/N$. Our main achievement is the proof of a local law for the individual elements of the resolvent - as well as various averaged versions - on these small scales. For this purpose, we rely on a dynamical approach (based on Dyson Brownian motion) leveraging the \emph{method of the characteristics} which was first introduced in \cite{lee2015edge} in the context of edge universality for deformed Hermitian matrices. The version of this method that we will use in this work is closer to the one introduced in \cite{bourgade2021extreme, huang2019rigidity}. More recently, this method has been extended and applied to a multitude of problems (see e.g. \cite{adhikari2020dyson, adhikari2023local, bourgade2206liouville, campbell2024spectral, cipolloni2023mesoscopic, cipolloni2024maximum, landon2022almost}). All these results prove local laws on scales $\eta\ge N^{-1+\xi}$, with the exceptions of \cite{cipolloni2024maximum} and \cite{huang2019rigidity}.  These works considered only the tracial local law on  logarithmic scales (with the latter still only working at scales much larger than the ones we require here) and were restricted to the spectral bulk. 

The main technical novelty of this work is to prove both an averaged and an isotropic local law on logarithmic scales $\eta\sim \log N/N$, uniformly throughout the spectrum. In particular, the scaling at the spectral edge introduces significant complications: the method of characteristics involves both the $\eta$ and $z$ parameters changing as $t$ changes. As a result, the relevant spectral scaling changes drastically as $z$ passes through the boundary of the spectral bulk (the unit circle). 

In order to overcome these difficulties, we find it convenient to track the evolution of the resolvent along characteristics in three separate stages. In the first stage we show that local laws with a polynomial error can be improved to a logarithmic error down to scales $\eta \geq N^{\xi-1}$. In the second stage we then show that this logarithmic error can be propogated down to relatively large logarithmic scales, $\eta \geq \frac{ ( \log N)^C}{N}$ for some $C>0$. Finally, we show that the error propogates from these large logarithmic scale to the short scale $\eta = \frac{ \log N}{N}$. The key consideration here is that the characteristic approach naturally introduces an error that scales like $\log (\eta_f / \eta_i)$ for a characteristic starting at an initial scale $\eta_i$ and ending at a final scale $\eta_f$. In the final step it is necessary to work with scales that are differentiated only by powers of $\log N$ as $\log ( \log N)^C = C \log \log N \ll \log N$. 


\subsection{Notations and conventions} For integers $k\in \nn$  we use the notation $[k]:=\{1, 2,\dots, k\}$.
For positive quantities $f, g$ we write $f\lesssim g$ and $f\asymp g$ if $f\le C g$ or $cg\le f\le Cg$, respectively, for some $N$--independent constants $c, C > 0$ which depend only on the constants appearing in \eqref{eq:moments} in the definition of our matrix model below. We denote vectors by bold-faced lower case Roman letters ${\bm x}, {\bm y}  \in\cc^d$, for some $d\in \nn$, and their scalar product by
\[
\langle {\bm x}, {\bm y}\rangle:=\sum_{i=1}^d \overline{x_i}y_i.
\]
For any $d \times d$ matrix $A$ we use the notation $\langle A\rangle:= d^{-1}\mathrm{Tr}[A]$ to denote the normalized trace of $A$, and $A^\mathfrak{t}$ denotes the transpose of $A$. We denote the $d$--dimensional identity matrix by $I=I_d$. Furthermore, we define the $2\times 2$ block matrices
\begin{equation}
\label{eq:defE1E2}
E_1:=\left(\begin{matrix}
1 & 0 \\
0 & 0
\end{matrix}\right),
\qquad\quad
E_2:=\left(\begin{matrix}
0 & 0 \\
0 & 1
\end{matrix}\right).
\end{equation}
We also use the notation
\[
\tilde{\sum}_{ij}:=\sum_{(i,j)\in \{(1,2),(2,1)\}},
\]
to denote sums over matrices $E_1,E_2$. For a vector ${\bf v}\in\cc^N$, with entries $v_i\in\cc$, and any $p\in\nn$ we define its $\ell^p$-norm by
\[
\lVert{\bf v}\rVert_p:=\left(\sum_{i=1}^N |v_i|^p\right)^{1/p},
\]
and $\lVert{\bf v}\rVert_\infty:=\max_i|v_i|$. For a matrix $A$ and vectors ${\bf x}$ and ${\bf y}$ we denote $A_{{\bf x}{\bf y}} := \langle {\bf x} , A {\bf y } \rangle$. 


We will use the concept of “with overwhelming probability” meaning that for any fixed $D > 0$ the probability of the event is bigger than $1-N^{-D}$ if $ N\geq N_0(D)$, with $N_0(D)$ possibly depending on the constants appearing in  \eqref{eq:moments} of the definition of our model, Definition \ref{def:model} below. Moreover, we use the convention that $\xi>0$ denotes an arbitrary small constant which is independent of $N$.


For real-valued martingales $M_t, N_t$, we denote the covariation process by $\d [ M_t, N_t]$. For complex valued martingales $M_t = X_t+ \i Y_t, N_t = P_t + \i Q_t$ the covariation process is defined by, $\d [ M_t, N_t] := \d [X_t, P_t] - \d [Y_t, Q_t] + \i ( \d [Y_t, P_t] + \d [ X_t, Q_t])$. The total variation process of a real-valued martingale is denoted by $[M_t] := \d [ M_t, M_t]$.

\medskip

\noindent\textbf{Acknowledgments.} The work of G.C. is partially supported by the MUR Excellence Department Project MatMod@TOV awarded to the Department of Mathematics, University of Rome Tor Vergata, CUP E83C18000100006. The work of B.L. is supported by NSERC and a Connaught New Researcher Award. 

\section{Main results}
\label{sec:mainres}

We consider $N\times N$ matrices $X$ with real or complex i.i.d. entries:
\begin{definition} \label{def:model} An \emph{i.i.d.\ matrix} is an $N \times N$ matrix $X$ whose entries are all independent, identically distributed (i.i.d.) random variables, $X_{ab} \stackrel{\d}{=} N^{-1/2} \chi$. On $\chi$ we assume that $\ee[\chi] = 0$ and $\ee[ |\chi|^2]=1$. We will consider two classes of i.i.d.\ matrices, real i.i.d.\ matrices and complex i.i.d.\ matrices. In the real case $\chi \in \rr$ and in the complex case $\chi \in \cc$ and we further assume that $\ee[ \chi^2]=0$. Additionally, we assume that for any $p \in \mathbb{N}$ there exists a $C_p >0$ so that,
\begin{equation}
\label{eq:moments}
\ee \left[ |\chi|^p \right] \leq C_p.
\end{equation}
Throughout, we will use the parameter $\beta$ to unify formulas that hold in the real and complex cases. Specifically, in the real case $\beta=1$ and in the complex case $\beta=2$. 
\end{definition}


For $\sigma\in \mathrm{Spec}(X)$ we denote its left and right eigenspaces by ${\bf L}_\sigma$ and ${\bf R}_\sigma$, and left and right eigenvectors in these spaces by ${\bf l}, {\bf r}$, respectively. Our main result is a sharp overwhelming probability bound on the size of these eigenvectors:
\begin{thm}
\label{theo:mainres}
Let $X$ be an i.i.d. matrix satisfying Definition~\ref{def:model}. For any $D>0$ there exists a constant $C= C(D) >0$ so that, for any unit vector ${\bf x} \in \cc^N$ we have,
\beq
\label{eq:highprobbound}
\pp\left( \sup_{\sigma \in \mathrm{Spec}(X)} \sup_{ {\bf l} \in {\bf L}_\sigma , {\bf r } \in {\bf R}_\sigma} \frac{|\langle{\bf x},{\bf l}\rangle|}{\lVert {\bf l}\rVert_2 }+\frac{|\langle {\bf x}, {\bf r}\rangle|}{\lVert {\bf r}\rVert_2 }\ge C\sqrt{\frac{\log N}{N}}\right)\le N^{-D}.
\eeq
\end{thm}

By sharpness, we mean that a decay in probability of the form $N^{-D}$ is not expected to hold if one replaces $\sqrt{ \log N}$ by any function that grows slower in $N$. In fact, the supremum in the probability is expected to be of order $\sqrt{ \frac{ \log N}{N}}$ with high probability. In this sense, the bound in \eqref{eq:highprobbound} is optimal in terms of the $N$--dependence. However, we do not attempt to quantify the relationship between $D$ and the constant $C$, and we do not expect that the simple method used here can access the smallest constant $C$ so that the probability on the left-hand side still tends to zero. 

Previously to \eqref{eq:highprobbound}, the best bound in term of $N$-dependence was $(\log N)^{9/2}/\sqrt{N}$ for matrices with sub--Gaussian entries, and $(\log N)^\gamma/\sqrt{N}$, for some large $\gamma>0$, for more general matrices \cite{MR3405592}.




\begin{rem}[More general ensembles]
We point out that we stated \eqref{eq:highprobbound} only for matrices with i.i.d. entries  to keep the presentation simple and short. However, our proof also applies for general matrices $X+A$, where $X$ is an i.i.d. matrix and  $A$ belongs to a large class of deterministic deformations (see e.g. \cite[Section 2.1]{cipolloni2024optimal} to study only bulk eigenvectors and \cite[Assumption 2.8]{campbell2024spectral} to include the edge). 
\end{rem}

\begin{rem}[Hermitian matrices]
In Theorem~\ref{theo:mainres} we focused on left/right eigenvectors of non-Hermitian matrices. However, it is easy to see that our proof applies with minor changes to eigenvectors of general Hermitian random matrices, giving an optimal overwhelming probability bound for eigenvectors of these matrices as well. For example, for Wigner matrices we could recover \cite[Theorem 1.2]{benigni2022optimal}, which was proven via a different method (the eigenvector moment flow). However, our simple approach does not allow us to compute the precise constant as in \cite[Theorems 1.3-1.4]{benigni2022optimal}.
\end{rem}

\section{Proof of Theorem~\ref{theo:mainres} by characteristics}
\label{sec:char}

\subsection{Hermitization and Matrix Dyson Equation}

A standard and useful tool for studying the spectrum of the non-Hermitian matrix $X$ is its Hermitization, which is defined by
\beq
\label{eq:herm}
H^z=H^z(X):=\left(\begin{matrix}
0 & X-z \\
(X-z)^* & 0
\end{matrix}\right), \qquad\quad z\in\cc.
\eeq
The eigenvalues of the matrix $H^z$ are related to those of $X$ by
\beq
z\in\mathrm{Spec}(X) \Longleftrightarrow 0\in\mathrm{Spec}(H^z).
\eeq
Note that due to its $2\times 2$ block structure, the spectrum of $H^z$ is symmetric with respect to zero. We denote its non--negative eigenvalues by $\{\lambda_i^z\}_{i\in [N]}$ (labeled in increasing order) and denote $\lambda_{-i}^z :=-\lambda_i^z$.

The main tool we will use for studying $H^z$ will be its resolvent, defined by 
\beq
\label{eq:defres}
G^z(\i\eta):=(H^z-\i\eta)^{-1}, \qquad\quad \eta>0.
\eeq
In fact, it is enough to give an upper bound for $\Im [G^z]$ that is uniform in $|z|$,  
since this will directly imply an upper bound for the non--Hermitian eigenvectors ${\bf l}, {\bf r}$ by the spectral theorem (see \eqref{eq:impbound} and the paragraph below it).

For large $N$  the resolvent $G^z$ is well approximated by the following matrix $M^z(\i\eta)$. To define it, we introduce $m^z ( \i \eta)$ as the unique solution to the Matrix Dyson Equation,
\beq
\label{eq:mde}
-\frac{1}{m^z(\i\eta)}=\i\eta+m^z(\i\eta)-\frac{|z|^2}{\i\eta+m^z(\i\eta)}, \qquad\quad \eta\Im m^z(\i\eta)>0.
\eeq
Then, given the solution of \eqref{eq:mde}, we define the $2N \times 2N$ matrix $M^z( \i \eta)$ by, (each block below is an $N \times N$ matrix)
\beq
\label{eq:mMrel}
M^z(\i\eta)=\left(\begin{matrix}
m^z(\i\eta)  & -zu^z(\i\eta) \\
-\overline{z}u^z(\i\eta) & m^z(\i\eta)
\end{matrix}\right), \qquad\quad u^z(\i\eta):=\frac{m^z(\i\eta)}{\i\eta+m^z(\i\eta)}.
\eeq
Theorem~\ref{theo:isolaw} below quantifies the above statement that $G^z$ is well approximated by $M^z$. Note that from \eqref{eq:mde} it follows that $m^z(\i \eta)$ is purely imaginary, 
and consequently that $u^z(\i\eta)$ is real. 
It will be useful to denote,
\beq \label{eqn:rho-asymp}
\rho^z ( \i \eta) := \Im[ m^z ( \i \eta)] \asymp \begin{cases} \eta^{1/3} + (1-|z|)^{1/2}, & |z| \leq 1 \\
\frac{\eta}{|z|-1+\eta^{2/3}} , & 10 > |z| \geq 1 \end{cases}
\eeq
with the asymptotics holding for all $0 < \eta < 1$ (the asymptotics follow from, e.g., \cite[Eq. (3.17)]{cipolloni2023universality}; see also \cite{alt2021spectral}).

\subsection{Local laws on short scales}

The main technical input to prove Theorem~\ref{theo:mainres} are the following local laws on short scales. The following should be compared with, e.g., \cite[Theorem 3.1]{cipolloni2024precise} where the logarithmic factors in the estimates are instead polynomial $N^\eps$ factors. In order to state them let us introduce the domain, for constants $C \geq 1$ and $0 < \xi < \frac{1}{100}$,
\beq
\D_{C, \xi} := \{ (z, \eta ) \in \cc \times (0, 1) : N \eta \rho^z ( \i \eta ) \geq 100 C^2 \log N , \eta / \rho^z ( \i \eta ) \leq N^{-\xi}, |z| \leq 1 + N^{-\xi} \}.
\eeq
\begin{thm}
\label{theo:isolaw}
Let $H^z$ and $G^z ( \i \eta) $ be as  in \eqref{eq:herm} and \eqref{eq:defres}. For any $D>0$ there is a constant $C= C(D) >0$ so that the following holds. Fix any  sufficiently small $\xi>0$, 
unit vectors ${\bf x}, {\bf y}\in\cc^{2N}$, and matrix $B\in\cc^{2N \times 2N}$. Then, 
\beq
\label{eq:avebound}
\pp\left[ \bigcap_{ (z, \eta) \in \D_{C, \xi}} \left\{ \big|\langle \big(G^z(\i\eta)-M^z(\i\eta)\big)B\rangle\big|\leq  \frac{C\lVert B\rVert\log N}{N\eta} \right\} \right] \geq 1 -  N^{-D},
\eeq
and
\beq
\label{eq:isobound}
\pp\left[ \bigcap_{ (z, \eta) \in \D_{C, \xi}} \left\{ \left|\langle {\bf x},\big(G^z(\i\eta)-M^z(\i\eta)\big){\bf y}\rangle\right|\leq C\sqrt{\log N}\sqrt{\frac{\rho^z ( \i \eta )}{N\eta}} \right\} \right] \geq 1 - N^{-D} ,
\eeq
for all $N$ large enough depending on $\xi$ and $D$.
\end{thm}

Strictly speaking to prove Theorem~\ref{theo:mainres} we only need \eqref{eq:isobound}. However, the proof of \eqref{eq:isobound} also requires control of the averaged quantities appearing on the LHS of \eqref{eq:avebound}, for the case $B=I$. Proving the estimate for general $B$ is not much harder and so is included above. 

The quantity $\eta \rho^z  ( \i \eta)$ plays an important role in our proofs and so we require the following simple lemma relating its asymptotics to other quantities.

\bel \label{lem:eta-rho}
Let $z$ satisfy $|z| \leq 10$. Let $1>\eta >0$ satisfy $ \eta \rho^z ( \i \eta) = A$ for $ 0 < A < 1$. Then,
\beq \label{eqn:char-final}
\eta \asymp \begin{cases}  A \frac{1}{ (1-|z|)^{1/2} +A^{1/4}} , & |z| \leq 1 \\  A^{1/2} ( \sqrt{ |z|-1} + A^{1/4} ) , & |z| \geq 1  \end{cases} , \qquad \frac{ \eta}{\rho^z ( \i \eta)} \asymp \begin{cases}  A \frac{1}{ (1-|z|) + A^{1/2}} & |z| \leq 1 \\  |z|-1+ A^{1/2}, & |z| \geq 1 \end{cases} ,
\eeq
and
\beq \label{eqn:rho-char-asymp}
\rho^z ( \i \eta) \asymp \begin{cases} A^{1/4} + \sqrt{1-|z|} , & |z| \leq 1 \\ \frac{A^{1/2}}{\sqrt{ |z|-1} + A^{1/4}} & |z| >1 \end{cases}.
\eeq
The function $\eta \to \eta \rho^z ( \i \eta)$ is strictly monotonically increasing and there is a constant $A_* >0$ so that for every $|z| \leq 10$ and $0 < A < A_*$ the equation $ \eta \rho^z ( \i \eta) = A$ has a unique solution for some $0 < \eta < 1$.
\eel
\proof One first establishes the asymptotics for $\eta$ using \eqref{eqn:rho-asymp}, and then substitutes the asymptotics for $\eta$ back into \eqref{eqn:rho-asymp} to find the asymptotics for the other two quantities. 

The monotonicity of $\eta \to \eta \rho^z ( \i \eta)$ follows from the fact that $m^z$ is the Stieltjes transform of a probability measure on $\rr$. The remaining statements follow from \eqref{eqn:rho-asymp}. 
\qed 

\

With the above we may now prove our main result.




\

\noindent{\bf Proof of Theorem~\ref{theo:mainres}.} Let $D>0$ and unit vector ${\bf x_1} \in \cc^N$ be given. Let $C=C(D+1)$ be the constant in Theorem \ref{theo:isolaw} and $\xi = \frac{1}{100}$. Let $\E$ be the event that the spectrum of $X$ lies within the disc of radius $1 + N^{-1/4}$ centered at the origin as well as the event on the LHS of \eqref{eq:isobound} for the domain $\D_{C, \xi}$ and the choice ${\bf x} = {\bf y } = ( 0 , {\bf x_1}^\mathfrak{t})^\mathfrak{t}$. Then by \eqref{eq:isobound} and \cite[Theorem 2.1]{alt2021spectral} we have that $\pp[ \E] \geq 1 - N^{-D}$.  We claim that the event in the probability on the LHS of \eqref{eq:highprobbound} holds on $\E$. The remainder of the proof is a deterministic argument verifying this statement and so we will assume that $\E$ holds. 

Let $\sigma\in \mathrm{Spec}(X)$, with left/right eigenvectors denoted by ${\bf l}$ and ${\bf r}$, respectively. Consider the vector
\beq
{\bf w}=\left(\begin{matrix} 0 \\ {\bf r}\end{matrix}\right)\in \cc^{2N} .
\eeq
Note that $H^\sigma {\bf w} = 0$, and so we may extend ${\bf w}/\lVert {\bf w} \rVert_2$ to a basis ${\bf w}/\lVert {\bf w}$ $ \rVert_2,$$ {\bf w}_2, \dots,$$ {\bf w}_N,$$ {\bf w}_{-1}, \dots$, $ {\bf w}_{-N}$ of $\cc^{2N}$ consisting of eigenvectors of $H^\sigma$, with corresponding eigenvalues $\lambda_1^\sigma=0,$ $ \lambda_2^\sigma, \dots,$ $ \lambda_N^\sigma, $ $\lambda_{-1}^\sigma, \dots, $ $  \lambda_{-N}^\sigma$. 
Then, denoting ${\bf x}=\left(\begin{matrix} 0 & {\bf x}^\mathfrak{t}_1 \end{matrix}\right)^\mathfrak{t}$, we have by the spectral theorem, for any $\eta>0$,
\beq
\label{eq:impbound}
\frac{|\langle{\bf x_1}, {\bf r}\rangle|^2}{\lVert {\bf r} \rVert_2^2}\leq \eta\left(\frac{|\langle{\bf x}, {\bf w}\rangle|^2}{\eta\lVert {\bf w} \rVert_2^2}+\eta\sum_{i \neq 1} \frac{|\langle {\bf x}, {\bf w}_i\rangle|^2}{(\lambda_i^\sigma)^2+\eta^2}\right)=\eta \langle {\bf x}, \Im G^\sigma(\i\eta){\bf x}\rangle.
\eeq
Choose $\eta=\eta(\sigma)$ so that $N\eta\rho^\sigma(\i \eta) =100C^2\log N$. Note that by Lemma  \ref{lem:eta-rho} such an $\eta$ exists, and that by \eqref{eqn:char-final} and the fact that $| \sigma| \leq 1 + N^{-1/4}$ we have that $\eta/ \rho^\sigma ( \i \eta) \leq N^{-1/5}$. Therefore,  $(\sigma, \eta) \in D_{C, \xi}$. By the fact that the estimates of \eqref{eq:isobound} hold on $\E$, we obtain the desired bound for $\langle{\bf x_1}, {\bf r}\rangle/\lVert {\bf r}\rVert_2$ by using \eqref{eq:isobound} to bound the RHS of \eqref{eq:impbound}.  

We conclude the proof by noticing that the same argument applies verbatim to the vector ${\bf w}=\left(\begin{matrix}{\bf l} & 0\end{matrix}\right)^\mathfrak{t}$.
\qed

\subsubsection{Dynamical approach to Theorem \ref{theo:isolaw} and characteristics}

The remainder of the paper is devoted to proving Theorem~\ref{theo:isolaw}. To prove \eqref{eq:isobound} we need to also prove an averaged version of this local law; for this reason, these two local laws will be proven in tandem. The proof of these local laws follows by a dynamical argument, i.e. we first prove these local laws for matrices with an (almost) order one Gaussian component using the \emph{method of characteristics} (see \cite{cipolloni2024maximum} for a similar proof in the averaged case) and then we show that this Gaussian component can be removed by a standard Green's function comparison argument.

Consider the Ornstein--Uhlenbeck flow
\beq
\label{eq:OU}
\d X_t=-\frac{1}{2}X_t \d t+\frac{\d B_t}{\sqrt{N}}, \qquad\quad X_0=X,
\eeq
with $X$ an i.i.d. matrix as defined in Definition~\ref{def:model}. Here $B_t$ is an $N\times N$ matrix of real or complex standard Brownian motions, for the cases $\beta=1$ and $\beta=2$, respectively. Let $H_t^z:=H^z(X_t)$ be the Hermitization defined as in \eqref{eq:herm} with $X$ replaced by $X_t$. In the remainder of the paper we denote the resolvent of $H_t^z$ by $G_t^z(\i\eta):=(H_t^z-\i\eta)^{-1}$, with $\eta>0$. Note that along the flow \eqref{eq:OU} the first and second moment of the entries of $X_t$, and so of $H_t^z$, are preserved along the whole flow. For this reason $M^z ( \i \eta)$ will remain a good approximation of the resolvent $G_t^z ( \i \eta)$. 
Our goal is now to study the evolution of $G_t^z(\i\eta)$ along the flow \eqref{eq:OU}. By the It\^{o} lemma we have, 
\begin{equation}
\label{eq:flownochar}
\begin{split}
\d \langle G^z_t(w)\rangle&=\d N^z_t(w)+\frac{1}{2}\langle G^z_t(w)\rangle \d t+\frac{1}{2}\langle (Z+w)G^z_t(w)^2\rangle \d t+2\tilde{\sum}_{ij}\langle G^z_t(w)E_i\rangle\langle G^z_t(w)^2E_j\rangle \d t \\
&\quad+ \frac{\1_{ \{ \beta=1\}}}{N} \tilde{\sum}_{ij} \langle G^z_t(w)^2 E_i G^z_t(w)^\mft E_j \rangle \d t
\end{split}
\end{equation}
and
\begin{equation}
\label{eq:flownochariso}
\begin{split}
&\d (G^z_t(w))_{\bf x y}=\d \widetilde{N}^z_t(w)+\frac{1}{2}(G^z_t(w))_{\bf x y} \d t+\frac{1}{2} (G^z_t(w)(Z+w)G^z_t(w))_{\bf x y} \d t\\
+ & 2\tilde{\sum}_{ij}\langle G^z_t(w)E_i\rangle (G^z_t(w)E_j G^z_t(w))_{\bf x y} \d t + \frac{\1_{ \{ \beta=1\}}}{N} \tilde{\sum}_{ij} (G^z_t(w) E_i G^z_t(w)^\mft E_jG^z_t(w))_{\bf x y} \d t.
\end{split}
\end{equation}
Here $\mft$ denotes the transpose of a matrix and $w=\i\eta$. The terms $\d N_t^z$ and $\d \widetilde{N}_t^z$ are Martingale increments defined in \eqref{eq:defbigBM} below (we need to introduce an additional equation so will introduce the Martingale increments all at once later). Note that in both these flows there is an additional term appearing in the real case $\beta=1$. 

The large $N$ limits of the flows \eqref{eq:flownochar} and \eqref{eq:flownochariso} are complex Burgers-type equations which may be studied via the method of characteristics. These characteristics are functions $ t \to (\eta_t, z_t) \in \rr \times \cc$, and are defined as the unique solution to the equations
\beq
\label{eq:defchar}
\partial_t \eta_t=-\Im m^{z_t}(\i\eta_t)-\frac{\eta_t}{2}, \qquad\quad \partial_t z_t=-\frac{z_t}{2},
\eeq
with initial conditions $(\eta_0, z_0)$ such that $\eta_0>0$ and $|z_0|\leq 10$.  A characteristic is defined up to some maximal time of existence $T_*$ when $\eta_{T_*} = 0$. Note that $m^{z_t}(\i\eta_t)$ depends on time only through $(\eta_t,z_t)$, that the characteristic for $\eta_t$ depends on $z_t$ through $m^{z_t}(\i\eta_t)$, and that
\beq
\label{eq:relM}
\partial_t M^{z_t}(\i\eta_t)=\frac{M^{z_t}(\i\eta_t)}{2}.
\eeq
In particular for $ 0 \leq t \leq T_* \wedge 1$ we see that,
\beq \label{eqn:rho-constant}
\rho_t \asymp \rho_0 .
\eeq
Furthermore, we notice that by \eqref{eq:defchar} and \eqref{eq:relM} for $s \leq t $ and $\rho_s:=\langle \Im M^{z_s}(\i\eta_s)\rangle$, we have 
\beq \label{eqn:del-t-ratio}
\del_t \left( \frac{\eta_t}{\rho_t} \right) = - 1 - \frac{\eta_t}{\rho_t} , \qquad \frac{\eta_t}{\rho_t}=\e^{s-t}\frac{\eta_s}{\rho_s}-(1-\e^{s-t}) .
\eeq
We require the following elementary lemma about the behavior of the characteristics.


\

\bel \label{lem:chars} 
Let $0.01 > \xi > 0$ and let $(\eta, z)$ be given. Define $A:= \eta \rho^z ( \i \eta)$ and assume that $|z| \leq 1 + N^{-10\xi}$ and $N^{-1} \leq A \leq N^{-20\xi}$. Let $T = N^{-\xi}$. Then there is a unique characteristic $( \eta_t, z_t)_{t \in [0, T]}$ so that $(\eta_T , z_T) = ( \eta, z)$ and the following holds.
\begin{enumerate}[label=\normalfont(\roman*)]
    \item We have $\eta_0 / \rho_0 \asymp N^{-\xi}$.
    \item Define the times $S_1 := T - \frac{N^{5\xi}}{N \rho_T^2}$ and $S_2 := T - \frac{ ( \log N)^3}{N \rho_T^2}$. Then $S_1 < S_2 < T$ and $T-S_1 \lesssim N^{-5 \xi}$.
    We have for $0 < s < S_1$ that,
    \beq \label{eqn:S1-char-bd}
N \eta_{s} \rho_s \gtrsim  N^{5 \xi},
    \eeq
    and for $0 < s < S_2$ we have,
    \beq
N \eta_s \rho_s \gtrsim  ( \log N)^3 .
    \eeq
    \item We have, 
    \beq  \label{eqn:chars-1}
\frac{\eta_0}{\eta_T} \gtrsim  N^{9 \xi}, \qquad \frac{\eta_{S_2}}{\eta_T} \lesssim  ( \log N)^3,
 \qquad \frac{ \eta_0}{\eta_{S_1}} \gtrsim c N^{4 \xi} .
    \eeq
\end{enumerate}
\eel
\proof The existence of a unique characteristic passing through $(\eta, z)$ follows by running \eqref{eq:defchar} backwards in time. The asymptotics in \eqref{eqn:char-final} together with the assumption $|z| -1 \leq N^{-10\xi}$ implies that $\eta / \rho ( \i \eta) \leq \O ( N^{-10\xi})$. Therefore the first item follows from \eqref{eqn:del-t-ratio}. 

From \eqref{eqn:rho-char-asymp} we see that $\rho_T \geq c N^{5 \xi-1/2}$ for some $c>0$ and so $\frac{N^{5\xi}}{N \rho_T^2} = \O ( N^{-5\xi} )$ showing that $T-S_1 = \O ( N^{-5 \xi} )$. From \eqref{eqn:rho-constant}--\eqref{eqn:del-t-ratio} one sees that for $0 < t < T$,
\beq \label{eqn:etarho-t}
\eta_t \rho_t \asymp \eta_T \rho_T + (T-t) \rho_T^2,
\eeq
and so the remaining estimates of the second item and the first two estimates of the third item follow easily. For the last estimate of \eqref{eqn:chars-1} we have $\eta_0 \gtrsim \rho_T N^{-\xi}$ and
\begin{align}
\eta_{S_1} \asymp \eta_T +\frac{ N^{5 \xi}}{N \rho_T} \lesssim  N^{-5 \xi} \rho_T,
\end{align}
and so the last estimate of \eqref{eqn:chars-1} follows. \qed

\

We will now study the evolution of the resolvent along the characteristics \eqref{eq:defchar}. In addition to the evolution of $\langle G_t^z (w) \rangle$,  we will - for technical reasons - also study the evolution of $\langle G_t^z(w)^2\rangle$. To lighten the notation we will use in the remainder of the paper the short-hand notation 
\beq
G_t:=G_t^{z_t}(\i\eta_t), \quad M_t:=M^{z_t}(\i\eta_t), \quad N_t:=N_t^{z_t}(\i\eta_t), \quad \widetilde{N}_t:=\widetilde{N}_t^{z_t}(\i\eta_t).
\eeq
Note in particular these quantities are evaluated at $z_t$ at time $t$. We also denote,
\beq
M_t' := \frac{\d}{\d w} M^{z_t} ( w ) \bigg\vert_{w=\i \eta_t} .
\eeq
For the above quantity we have the estimate,
\beq \label{eqn:M'-bd-a}
|\langle M_u' \rangle| \leq \frac{ \rho_u}{\eta_u}
\eeq
which follows since $\langle M^z (w) \rangle$ is the Stieljes transform of a probability measure. 
We have the equations, 
\begin{equation}
\label{eq:flownocharchar}
\begin{split}
\d \langle G_t-M_t\rangle &=\d N_t+\Phi(t)\langle G_t-M_t\rangle \d t+ \langle G_t-M_t\rangle\langle G_t^2-M_t'\rangle \d t+ \frac{\1_{ \{ \beta=1\}}}{N} \tilde{\sum}_{ij} \langle G_t^2 E_i G_t^\mft E_j \rangle \d t, \\
\d \langle G_t^2-M_t'\rangle &=\d \widehat{N}_t+2\Phi(t)\langle G_t^2-M_t'\rangle \d t+\langle G_t^2-M_t'\rangle^2 \d t+2\langle G_t-m_t\rangle\langle G_t^3\rangle\d t \\
&\quad+ \frac{\1_{ \{ \beta=1\}}}{N} \tilde{\sum}_{ij} \big(\langle G_t^3 E_i G_t^\mft E_j \rangle+\langle G_t^2 E_i (G_t^2)^\mft E_j \rangle\big) \d t
\end{split}
\end{equation}
and
\begin{equation}
\label{eq:flownocharisochar}
\d (G_t-M_t)_{\bf x y}=\d \widetilde{N}_t+\frac{1}{2}(G_t-M_t)_{\bf x y} \d t+\langle G_t-M_t\rangle (G_t^2)_{\bf x y} \d t + \frac{\1_{ \{ \beta=1\}}}{N} \tilde{\sum}_{ij} (G_t E_i G_t^\mft E_jG_t)_{\bf x y} \d t 
\end{equation}
where we used \eqref{eq:relM} and that $2\langle G^{z_t}_t(w)E_i\rangle=\langle G^{z_t}_t(w)\rangle$. 
Here, 
\beq \label{eqn:Phi-def}
\Phi(t):=\frac{1}{2}+\langle M_t'\rangle,
\eeq
and the stochastic terms are defined by
\begin{equation}
\label{eq:defbigBM}
\begin{split}
\d N^z_t(w):&=-\frac{1}{\sqrt{N}}\langle (G^z_t(w))^2\d\mathfrak{B}_t \rangle, \qquad \d \widehat{N}^z_t(w):=-\frac{1}{\sqrt{N}}\langle (G^z_t(w))^3\d\mathfrak{B}_t \rangle \\
\d \widetilde{N}^z_t:&=-\frac{1}{\sqrt{N}}(G^z_t(w)\d\mathfrak{B}_tG^z_t(w))_{\bf x y} , \quad \mathfrak{B}_t := \left( \begin{matrix} 0 & B_t \\ B_t^* & 0 \end{matrix}\right).
\end{split}
\end{equation}


The following proposition will be the main technical input for the proof of the local laws in Theorem~\ref{theo:isolaw}.
\begin{prop}
\label{pro:zig} Fix any $D>0$. There is a $K>0$ (depending on $D>0$) so that the following holds. 
Fix any sufficiently small $\xi>0$ and denote $T = N^{- \xi}$. Let $( \eta_t, z_t)_{0 < t < T}$ denote a characteristic such that $K N^{-1} \log N  \leq \eta_T \rho_T \leq N^{-100 \xi}$ and $|z_T| \leq 1 + N^{-100\xi}$. 
Let $B\in\cc^{2N\times 2N}$ be a deterministic matrix, and let ${\bf x}, {\bf y}\in\cc^{2N}$ be unit vectors. Then, for all $t \in [ T- N^{5 \xi}/ ( N \rho_T^2), T]$ we have 
\beq
\label{eq:avegtlaw}
\pp\left(\big|\langle B(G_t^{z_t}(\i\eta_t)-M^{z_t}(\i\eta_t))\rangle\big|\geq \frac{K\lVert B\rVert \log N}{N\eta_t}\right)\leq N^{-D},
\eeq
and
\beq
\label{eq:isogtlaw}
\pp\left(\left|\langle {\bf x},\big(G_t^{z_t}(\i\eta_t)-M^{z_t}(\i\eta_t)\big){\bf y}\rangle\right|\geq K \sqrt{\log N}\sqrt{\frac{\rho_t}{N\eta_t}}\right)\leq N^{-D} .
\eeq
\end{prop}


The proof of Proposition~\ref{pro:zig} appears in the next section. We finally conclude the proof of Theorem~\ref{theo:isolaw}.

\

\noindent{\bf Proof of Theorem~\ref{theo:isolaw}.} We first recall here a basic argument based on continuity which says that it suffices to prove  \eqref{eq:avebound}--\eqref{eq:isobound} for a single $z$ and $\eta$; i.e., without the intersection over all $(z, \eta) \in \D_{C, \xi}$. To be precise, uniformly in the parameters $(z, \eta)$ with $\eta \geq N^{-1}$ we have that the matrix entries of $G^z ( \i \eta)$ and $M^z ( \i \eta)$ are jointly Lipschitz in the parameters $z, \eta$ with constant bounded above by $N^{2}$. In particular, we can find a set $\mathcal{S}$ of at most $N^{10}$ points $(z, \eta) \subseteq \D_{C, \xi}$ such that for every $z$ and $\eta$ inside the probabilities in \eqref{eq:avebound}, \eqref{eq:isobound} we have
\beq \label{eqn:grid-est}
G_{ij}^{z} ( \i \eta)- M_{ij}^{z } ( \i \eta ) = G_{ij}^{w} ( \i \theta) - M_{ij}^{w } ( \i \theta ) + \O ( N^{-3})
\eeq
for some $(w, \theta) \in \mathcal{S}$ satisfying $|w-z| + | \theta - \eta| \leq N^{-5}$. Therefore, if we can prove \eqref{eq:avebound} and \eqref{eq:isobound} for fixed $z$ and $\eta$, the theorem will follow by a union bound and \eqref{eqn:grid-est}.

To prove the estimate for a fixed $z$ and $\eta$, we first note that Proposition~\ref{pro:zig} proves the bounds \eqref{eq:avebound}--\eqref{eq:isobound} for matrices with a Gaussian component of size $T=N^{-\xi}$; that is, for matrices of the form $X = (1- \e^{-T})^{1/2} G + \e^{-T/2} Y$ with $G$ a matrix of i.i.d. Gaussian random variables and $Y$ a general i.i.d. matrix satisfying Definition \ref{def:model}. What therefore remains to be shown is that if the estimate holds for this general class of so-called \emph{Gaussian divisible matrices}, then the same estimate holds for any i.i.d. matrix.  This can be deduced from a modification of what is known as the four moment method, or Green's function comparison theorem \cite[Theorem 15]{tao2011random} which is a standard approach in random matrix theory that allows one to compare spectral observables of two different random matrix ensembles. In the work \cite{cipolloni2024maximum} we carried out almost this exact argument for the observable on the left-hand side of \eqref{eq:avebound} in the case that $B = \1$. In Appendix~\ref{sec:conclproof} we detail what changes for the case of general $B$ as well as for the isotropic observables on the left-hand side of \eqref{eq:isobound}.  \qed



\





\section{Proof of Proposition~\ref{pro:zig}}

The proof of Proposition \ref{pro:zig} will be divided into three main parts. In the first part we prove the local law \eqref{eq:avegtlaw} for the case $B=\1$. This is the most delicate part of the proof. Further, the estimates obtained will be used in the proof of \eqref{eq:avegtlaw} for general $B$ as well as in \eqref{eq:isogtlaw}. Precisely, Proposition \ref{pro:zig} follows from Proposition \ref{prop:ll-id}, Lemma \ref{lem:B-flow} and Lemma \ref{lem:iso-flow} below.

Before starting with the proofs we recall that our characteristic $(\eta_t, z_t)_{0 < t < T }$ satisfies $K_0 N^{-1} \log N \leq \rho_T \eta_T \leq N^{-100\xi}$, $|z_T| \leq 1 + N^{-100\xi}$ and $T = N^{-\xi}$. From \eqref{eq:defchar} and \eqref{eq:relM} it follows that for $t < T$, $\eta_t \asymp \eta_0 + (T-t) \rho_T$. In particular, by applying \eqref{eqn:rho-asymp}, \eqref{eqn:char-final} and letting $A:= \eta_T \rho_T$  we have that,
\beq
\label{eq:linetarhno}
\eta_t \rho_t \asymp \eta_T \rho_T + (T-t) \times \begin{cases} \frac{A}{ |z_T|-1 +A^{1/2}}, & |z_T| >1 \\ A^{1/2} + (1-|z_T|) , & |z_T| \leq 1 \end{cases}.
\eeq
Additionally, from Lemma \ref{lem:chars} we note that $\eta_0 / \rho_0 \asymp N^{-\xi}$, and with
\beq
\label{eq:defs1s2}
S_1= T-\frac{N^{5\xi}}{N\rho_T^2}, \qquad\quad S_2= T -\frac{(\log N)^3}{N\rho_T^2},
\eeq
we have
$N\eta_t \rho_t \gtrsim N^{5 \xi}$, for $0 \leq t \leq S_1$, and $N\eta_t \rho_t \gtrsim (\log N)^3$, for $S_1\leq t \leq S_2$.

\bep \label{prop:ll-id}
For any $D>0$, if $K_0$ is large enough we have that for $t \in [S_1, T]$ that
\beq \label{eqn:prop-ll-av}
\pp\left[ |\langle G_t - M_t \rangle| \geq \frac{K_0 \log N}{N \eta_t} \right] \leq N^{-D}
\eeq
for $N$ large enough depending on $D, \xi$. 
\eep
The proof of Proposition \ref{prop:ll-id} is itself divided into three parts. We will analyze the evolution of $G_t$ in the three time intervals $[0, S_1]$, $[S_1, S_2]$ and $[S_2, T]$.  In the first time interval we show that we can improve the error term in the local laws for times $t \approx S_1$ from a polynomial error $\frac{N^\xi}{N \eta}$ to a logarithmic one $\frac{\log N}{N \eta}$.  Note that for such times the spectral parameter $\eta_t$ is relatively large (as $N\eta_t\rho_t\ge N^{5\xi}$ as seen above). In the second time interval of the proof we propagate this better error starting at time $t = S_1$ to $t=S_2$, when $\eta_t$ is smaller, i.e. $N\eta_t\rho_t\ge (\log N)^3$, but still not on the desired logarithmic scale. Finally, in the third part we propagate this better error from time $t=S_2$ to $t=T$ when we reach $N\eta_t\rho_t\ge C_1 \log N$. We point out that this proof follows a similar strategy to \cite[Lemma 3.4]{cipolloni2024maximum}, with the important difference that in the proof of \cite[Lemma 3.4]{cipolloni2024maximum} only the first two parts were needed, since we considered only the bulk of the spectrum $|z_T|<1$. Instead, we now consider also the edge of the spectrum leading to the additional difficulty of taking into account the smallness of $\rho_T$ when $|z_T|\approx 1$, which requires the third additional step.

Before beginning the proofs let us introduce some further notation. We define,
\beq
\label{eq:shortnot}
X_{1}^{\av}(t) :=\langle G_t-M_t\rangle, \qquad\quad X_{2}^{\av}(t):=\langle G_t^2-M_t'\rangle,
\eeq
as well as,
\beq \label{eqn:tstar-def}
t_* := \inf\{ t \geq 0 : |z_t| \leq 1 \},
\eeq
the first time the characteristic enters the unit disc (note that it is $0$ if the characteristic starts inside the disc). Note that by \eqref{eq:defchar} if $|z_0 | > 1$ we have
\beq \label{eqn:tstar-asymp}
t_* = 2 ( |z_0| -1) + \O ( (|z_0| -1 )^2 ).
\eeq
We also require the parameter,
\beq \label{eqn:kap-def}
\kappa_0:= (|z_0|-1)_+^{3/2} + \eta_0. 
\eeq

\bel \label{lem:av-ll-1} Let $D>0$ be given. There is a $K_1 =K_1 (D) >0$ so that with probability at least $1- N^{-D}$ for all $t \in [0, S_1]$ we have
\beq \label{eq:impsingG}
| \langle G_t - M_t \rangle | \leq  \frac{ K_1 \log N}{N \eta_t} +  \frac{ N^{\xi}}{N \eta_{t_* \wedge t}} \frac{ \eta_0}{\kappa_0}
\eeq
 and
\beq \label{eq:imp-G-1}
| \langle G_{S_1} - M_{S_1} \rangle |\leq  \frac{ K_1 \log N}{N \eta_{S_1}} , \qquad  | \langle G_{S_1}^2 - M'_{S_1}\rangle | \leq  \frac{K_1 ( \log N)^2}{N \eta_{S_1}^2}.
\eeq
\eel
\proof Define the stopping time,\beq
\label{eq:deftau1}
\tau_1:=\inf \left\{t\geq 0 : |X_{1}^{\av} (t)|\geq\frac{N^{\xi}}{N\eta_t},  \quad |X_{2}^{\av} (t) |\geq \frac{N^{2\xi}}{N\eta_t^2}\right\}\wedge S_1.
\eeq
Recall that $N\eta_t\rho_t\geq N^{5\xi}$ for $0\le t\le S_1$.  Then by \cite[Eq. (3.15)]{campbell2024spectral}\footnote{For the local law in \eqref{eq:kllaws} we refer to \cite{campbell2024spectral}; however, this bound could easily be proven from scratch following a fairly standard cumulant expansion, since we need it only for $\kappa_0\gtrsim N^{-\xi}$. The same remark also applies to the estimate of the initial condition in \eqref{eqn:inB} and \eqref{eq:iniso} below.} we have\footnote{We point out that \cite[Eq. (3.15)]{campbell2024spectral} proves \eqref{eq:kllaws} only for $|X_{1}^{\av} (0)|$. However, given this bound obtaining the desired bound for $|X_{2}^{\av}(0)|$ is immediate using Cauchy's integral formula as in \cite[Eq. (8.63)]{campbell2024spectral}.}
\beq
\label{eq:kllaws}
|X_{1}^{\av}(0)|\leq \frac{N^{\xi/2}}{N\kappa_0},  \qquad\quad |X_{2}^{\av} (0)|\leq \frac{N^{\xi/2}}{N\kappa_0^2},
\eeq
with overwhelming probability. This shows that $\tau_1>0$ with overwhelming probability. In fact one can check that $\tau_1 = S_1$ from \cite[Theorem 5.2]{alt2021spectral}  (and an  analogous estimate for $\langle G^2\rangle$), but we will not use this directly as the analog will not be true in the proof of Lemmas \ref{lem:av-ll-2} and \ref{lem:av-ll-3} below. 

Notice that for $t<\tau_1$ and $l\geq 2$, by Schwarz inequality and Ward identity, we have for $N$ sufficiently large,
\beq
\label{eq:aviso}
\left| \langle G_t^l\rangle \right| \leq \langle |G_t|^l\rangle\leq \frac{1}{\eta_t^{l-1}}\left(\rho_t+\frac{N^{\xi}}{N\eta_t}\right)\leq 2 \frac{\rho_t}{\eta_t^{l-1}}, 
\eeq
where we used \eqref{eqn:S1-char-bd}. 
Then, integrating \eqref{eq:flownocharchar} in time and using \eqref{eq:aviso} and \eqref{eqn:schwarz-ave}, we obtain
\begin{align}
\label{eq:flowgro}
\big|X_{1}^\av (t \wedge \tau_1 ) \big|&\leq \big|p_{0,t\wedge\tau_1}X_{1}^\av (0) \big|+\left|\int_0^{t\wedge\tau_1}p_{s,t\wedge\tau_1}\,\d N_s\right|+\int_0^{t\wedge\tau_1} p_{s,t\wedge\tau_1}\big|X_{1}^\av (s) X_{2}^\av (s) \big|\,\d s \notag \\
&\quad+4\int_0^{t\wedge\tau_1} p_{s,t\wedge\tau_1}\frac{\rho_s}{N\eta_s^2}\,\d s, 
\end{align}
and
\begin{align}\label{eq:flowgro-2}
\big|X_{2}^\av (t \wedge \tau_1 ) \big|&\leq \big|p_{0,t\wedge\tau_1}^2X_{2}^\av (0) \big|+\left|\int_0^{t\wedge\tau_1}p_{s,t\wedge\tau_1}^2\,\d \widehat{N}_s\right|+\int_0^{t\wedge\tau_1} p_{s,t\wedge\tau_1}^2\big| X_{2}^\av (s) \big|^2\,\d s \notag \\
&\quad+2\int_0^{t\wedge\tau_1} p_{s,t\wedge\tau_1}^2\left(\frac{\rho_s}{\eta_s^2}\big|X_{1}^\av (s) \big|+4\frac{\rho_s}{N\eta_s^3}\right)\,\d s, 
\end{align}
where (recall the definition of $\Phi$ in \eqref{eqn:Phi-def})
\beq
\label{eq:propag}
p_{s,t}:=\exp\left(\int_s^t\Phi(r)\,\d r\right)\leq 2\frac{\eta_s}{\eta_t}.
\eeq
This inequality follows by $\Phi(r)\le 1/2+\rho_t/\eta_t$ as a consequence of \eqref{eqn:M'-bd-a} above. 

For $N$ sufficiently large we have
\begin{align}
\label{eqn:aaa-2}
& \pp\left[\sup_{0<u<S_1}N\eta_u\left|\int_0^{u\wedge\tau_1}p_{r,u\wedge\tau_1}\,\d N_r\right| >  K \log N \nc\right] \leq N^{-D} , \notag\\
& \pp\left[ \sup_{0<u<S_1}N\eta_u^2\left|\int_0^{u\wedge\tau_1}p_{r,u\wedge\tau_1}^2\,\d \widehat{N}_r\right| >  K \log N \nc\right]\le N^{-D},
\end{align}
for any $D>0$ and for an $N$--dependent $K=K(D)\geq 100$ appropriately chosen. The proof of these two inequalities follows immediately by  Lemma~\ref{lem:stochastic}.

Using the definition of $\tau_1$ and \eqref{eq:propag}, we see that the last two terms on the RHS of \eqref{eq:flowgro} are bounded by
\begin{equation}
\label{eq:usbound}
\begin{split}
\int_0^{t\wedge\tau_1} p_{s,t\wedge\tau_1}\big|X_{1}^\av (s) X_{2}^\av  (s) \big|\,\d s&\leq 2\int_0^{t\wedge\tau_1} \frac{\eta_s}{\eta_{t \wedge \tau_1}}\frac{N^{\xi}}{N\eta_s}\frac{N^{2\xi}}{N\eta_s^2}\leq \frac{1}{N\eta_{t \wedge \tau_1}}, \\
\int_0^{t\wedge\tau_1}p_{s,t\wedge\tau_1}\frac{\rho_s}{N\eta_s^2}\,\d s&\le 2 \int_0^{t\wedge\tau_1}\frac{\rho_s}{N\eta_s\eta_{t\wedge\tau_1}}\,\d s\leq 3 \frac{\log N}{N\eta_{t\wedge\tau_1}},
\end{split}
\end{equation}
where in the first line we used \eqref{eqn:S1-char-bd}, and \eqref{eq:defchar} was used to estimate the time integrals in the last inequalities on both lines.  
To estimate the initial condition in \eqref{eq:flowgro} (i.e., the first term on the RHS) we will use the stronger bound
\beq
\label{eq:imppropb}
p_{s,t}\leq C \frac{\eta_{s\wedge t_*}}{\eta_{t\wedge t_*}},
\eeq
where $t_*$ is defined in \eqref{eqn:tstar-def}. 
The proof of this bound appears in Lemma \ref{lem:av-ll-1a} below. In particular, \eqref{eq:imppropb} and \eqref{eq:kllaws} imply that,
\beq
\label{eq:incondest}
\big|p_{0,t\wedge\tau_1}X_{1}^\av (0)\big| 
\leq \frac{C N^{\xi/2}}{N\eta_{t_*\wedge t\wedge\tau_1}}\frac{\eta_0}{\kappa_0} .
\eeq
 Using the estimates \eqref{eqn:aaa-2}, \eqref{eq:usbound}, \eqref{eq:incondest} in \eqref{eq:flowgro} 
 we conclude
\beq \label{eq:impsingG-a}
|X_{1}^{\av} (t \wedge \tau_1) |\leq \frac{2 K \log N}{N\eta_{t\wedge\tau_1}} + \frac{C N^{\xi/2}}{N\eta_{t_*\wedge t\wedge\tau_1}}\frac{\eta_0}{\kappa_0}.
\eeq
We now turn to estimating the RHS of \eqref{eq:flowgro-2}. We estimate its last two terms by
\begin{equation}
\label{eq:usbound2}
\begin{split}
\int_0^{t\wedge\tau_1} p_{s,t\wedge\tau_1}^2\big| X_{2}^\av (s) \big|^2\,\d s&\leq 2 \int_0^{t\wedge\tau_1}\frac{\eta_s^2}{\eta_{t\wedge\tau_1}^2}\frac{N^{4\xi}}{N^2\eta_s^4}\,\d s\leq \frac{1}{N\eta_{t \wedge \tau_1}^2} \\
\int_0^{t\wedge\tau_1} p_{s,t\wedge\tau_1}^2\left(\frac{\rho_s}{\eta_s^2}\big|X_{1}^\av (s) \big|+4\frac{\rho_s}{N\eta_s^3}\right)\,\d s &\leq \int_0^{t\wedge\tau_1}\frac{\eta_{s\wedge t_*}^2}{\eta_{t_*\wedge t\wedge\tau_1}^2}\left(\frac{CK \rho_s\log N}{N\eta_s^3}+ \frac{\rho_s}{\eta_s^2}\frac{C N^{\xi/2}}{N\eta_{t_*\wedge s}}\frac{\eta_0}{\kappa_0}\nc\right)\,\d s\\
&\le C K \frac{(\log N)^2}{N\eta_{t\wedge\tau_1}^2}+\frac{C N^{\xi/2} \log N }{N\eta_{t_*\wedge t\wedge\tau_1}^2}\frac{\eta_0}{\kappa_0}.
\end{split}
\end{equation}
In the first line we used \eqref{eqn:S1-char-bd}. In the second line we used \eqref{eq:impsingG-a} and \eqref{eq:imppropb}. 
For the initial condition (the first term on the RHS of \eqref{eq:flowgro-2}) we have
\beq \label{eqn:X2-init}
\big|p_{0,t\wedge\tau_1}^2 X_{1}^\av (0) \big|
\le  \frac{N^\xi}{N\eta_{t_*\wedge t\wedge\tau_1}^2}\left(\frac{\eta_0}{\kappa_0}\right)^2,
\eeq
where we used \eqref{eq:imppropb} and \eqref{eq:kllaws}.

Using \eqref{eqn:aaa-2}, \eqref{eq:usbound2} and \eqref{eqn:X2-init} in \eqref{eq:flowgro-2} we conclude that
\beq
\label{eq:impsingG2}
|X_{2}^{\av} (t \wedge \tau_1 )|\le \frac{C K (\log N)^2}{N\eta_{t\wedge\tau_1}^2} +\frac{N^\xi}{N\eta_{t_*\wedge t\wedge\tau_1}^2}\left(\frac{\eta_0}{\kappa_0}\right)^2+\frac{N^\xi }{N\eta_{t_*\wedge t\wedge\tau_1}^2}\frac{\eta_0}{\kappa_0},
\eeq
with probability at least $1- N^{-D}$. The bounds \eqref{eq:impsingG-a} and \eqref{eq:impsingG2} show that $\tau_1=S_1$ with probability at least $1 - 2N^{-D}$. Therefore, \eqref{eq:impsingG} holds by \eqref{eq:impsingG-a}.  Furthermore, choosing $t=S_1$ and using that (the proof of this bound is presented below in Lemma \ref{lem:av-ll-1a})
\beq
\label{eq:ei}
\frac{\eta_0}{\eta_{S_1\wedge t_*}\kappa_0}\lesssim \frac{N^{-2\xi}}{\eta_{S_1}}
\eeq
in \eqref{eq:impsingG-a} and \eqref{eq:impsingG2},
we obtain
\beq
\label{eq:newimpbs}
|X_{1}^{\av} (S_1)|\le\frac{3 K \log N}{N\eta_{S_1}}, \qquad\quad |X_{2}^{\av} (S_1)|\le\frac{C K(\log N)^2}{N\eta_{S_1}^2}.
\eeq
This completes the proof. \qed

\bel \label{lem:av-ll-1a}
In the notation of the above lemma we have, for some $C>0$,
\beq \label{eq:incondest-1}
p_{s, t} \leq C \frac{ \eta_{s \wedge t_*}}{\eta_{t \wedge t_*}}
\eeq
and
\beq \label{eq:ei-1}
\frac{ \eta_0}{ \eta_{S_1 \wedge t_*} \kappa_0} \leq C \frac{ N^{-2 \xi}}{\eta_{S_1}}
\eeq
\eel
\proof We start with the proof of the first inequality. We first have (see Appendix \ref{a:M'-bd})
\beq \label{eqn:M'-bd}
\big|\langle M_r'\rangle\big|\le \frac{1}{2\rho_r^2+\eta_r/\rho_r}\left(1+\mathcal{O}\left(\frac{\eta_r}{\rho_r}\right)\right).
\eeq
Then, we split the estimate of the exponential of the time integral in \eqref{eq:incondest-1} into two regimes: $r\in [s,t_*]$ and $r\in (t_*, t]$. Notice, that in certain cases one of these two regimes may be empty. 

We start with the regime $r\in (t_*,t]$. Here, $\rho_r \geq c \eta_r$ by \eqref{eqn:rho-asymp} and so $| \langle M_r' \rangle| \lesssim \rho_r^{-2}$. 
We estimate,
\beq
\exp\left(\int_{t_*}^t \left[\frac{1}{2}+\langle M_r'\rangle\, \d r\right]\right)\le 2 \exp\left(\int_{t_*}^t \frac{C}{\rho_r^2}\, \d r\right)\le 2 \exp\left( C \frac{t-t_*}{\rho_{t_*}^2}\right)\lesssim  1.
\eeq
Here in the last inequality we used $t-t_*\lesssim \eta_{t_*}/{\rho_{t_*}}$ by \eqref{eqn:del-t-ratio}, and that $\eta_{t_*}\lesssim \rho_{t_*}^3$ by \eqref{eqn:rho-asymp} as $|z_{t_*}| \leq 1$.

In the regime $r\in [s,t_*]$ (which is present only if $|z_0|>1$) we estimate using simply \eqref{eqn:M'-bd}
\beq
\exp\left(\int_s^{t_*} \left[\frac{1}{2}+\langle M_r'\rangle\, \d r\right]\right)\le 2 \exp\left(\int_s^{t_*} \frac{\rho_r}{\eta_r}\, \d r\right)\le2 \frac{\eta_s}{\eta_{t_*}}.
\eeq
This concludes the proof of \eqref{eq:incondest-1}. 

We now turn to the proof of \eqref{eq:ei-1}. We first consider the case that $t_* = 0$ so that $|z_0| \leq 1$. Then $\kappa_0 = \eta_0$ and the desired inequality follows from the fact that $\eta_0 \geq N^{2 \xi} \eta_{S_1}$, which follows from \eqref{eqn:chars-1}. 

In the remainder of the proof we therefore assume that $|z_0| > 1$. If $\kappa_0 \geq N^{2 \xi} \eta_0$, then \eqref{eq:ei-1} is immediate so we may further assume that $\kappa_0 \leq N^{2 \xi} \eta_0$. We first consider the case $t_* < S_1$. Then from $\rho_0 \asymp \frac{ \eta_0}{ |z_0|-1 + \eta_0^{2/3}}$ (see \eqref{eqn:rho-asymp}), the definition of $\kappa_0$ in \eqref{eqn:kap-def}, and the fact that $\eta_{t_*} \asymp \rho_{t_*}^3 \asymp \rho_0^3$ we have,
\beq
\frac{\eta_0}{\eta_{t_*}\kappa_0}\asymp \frac{1}{\rho_0^{3/2}\eta_0^{1/2}}\lesssim \frac{N^{2\xi}}{\eta_0}\lesssim\frac{N^{-2\xi}}{\eta_{S_1}}.
\eeq
We point out that in the penultimate inequality we used that in this regime  $\rho_0^{3/2}\gtrsim N^{-2\xi}\eta_0^{1/2}$ (as a consequence of $\kappa_0\asymp (\eta_0/\rho_0)^{3/2}$ and $\kappa_0\leq N^{2\xi}\eta_0$). 
The last inequality  used $\eta_0 \gtrsim N^{4 \xi} \eta_{S_1}$ from \eqref{eqn:chars-1}. 
This yields \eqref{eq:ei-1} in this case.

Finally, we consider the case $t_* \geq S_1$ and $\kappa_0 \leq N^{2 \xi} \eta_0$. Note that by Lemma \ref{lem:chars} we have that $S_1 \asymp T$ and so by \eqref{eqn:tstar-asymp} we get that $|z_0| -1 \asymp N^{-\xi}$ and so 
\beq
\kappa_0 \geq c N^{-3 \xi/2}.
\eeq
We claim that $\eta_0 \leq N^{-10 \xi}$ which together with the above inequality will complete the proof of \eqref{eq:ei-1}. To see that $\eta_0 \leq N^{-10 \xi}$ we first note that
\beq
\eta_0 \asymp N^{-\xi} \rho_T
\eeq
by Lemma \ref{lem:chars}. If $|z_T| >1$, since $\rho_T \eta_T \leq N^{-100\xi}$ by assumption and $\eta_T \geq c \rho_T^{3}$ by \eqref{eqn:rho-asymp} (assuming $|z_T| > 1$) we get $\rho_T \leq C N^{-25 \xi}$, which implies $\eta_0 \leq N^{-10\xi}$ as desired. 

On the other hand, if $|z_T| \leq 1$, by \eqref{eqn:rho-char-asymp}, we have
\begin{align}
\rho_T &\asymp \sqrt{ 1 -|z_T|} + (\eta_T \rho_T )^{1/4} \asymp \sqrt{ T-t_*} + ( \eta_T \rho_T)^{1/4} \notag\\
& \leq \sqrt{T-S_1} + N^{-25\xi} \asymp \frac{ N^{5 \xi/2}}{N^{1/2} \rho_T} + N^{-25 \xi} \leq N^{5 \xi/2} \sqrt{ \frac{\eta_T}{\rho_T}} + N^{-25 \xi}  \notag\\
& \leq \frac{ N^{5 \xi/2} N^{-50\xi}}{ \rho_T} + N^{-25\xi}
\end{align}
The last inequality on the second line follows from the assumption $\eta_T \rho_T \geq N^{-1}$. The third line follows from $\eta_T \leq N^{-100\xi} \rho_T^{-1}$. This final inequality implies that $\rho_T \leq N^{-20 \xi}$ and so the proof is completed. \qed

\

We now show that the bounds \eqref{eq:newimpbs},
which hold at time $S_1$,  can be propagated down to $N\eta_t\rho_t\ge (\log N)^3$, for $S_1\le t\le S_2$.

\bel \label{lem:av-ll-2}
For any $D>0$ there is a $K_2 = K_2(D)$ so that with probability at least $1  - N^{-D}$ we have,
\beq \label{eq:imprllaw}
|X_1^\av (t) | \leq \frac{ K_2 \log N}{N \eta_t} , \qquad | X_2^\av (t) | \leq \frac{ K_2 ( \log N)^2}{N \eta_t^2} ,
\eeq
for all $ t\in [S_1, S_2]$. 
\eel
\proof The proof of this is analogous to Lemma \ref{lem:av-ll-1} so we focus on explaining the differences. 
We replace the stopping time $\tau_1$ from \eqref{eq:deftau1}  with
\beq \label{eqn:tau-2-def}
\tau_2:=\inf \left\{t\ge S_1 : |X_{1}^{\av} (t) |=\frac{5K\log N}{N\eta_t}, \quad |X_{2}^{\av} (t) |=\frac{20K(\log N)^2}{N\eta_t^2}\right\}\wedge S_2.
\eeq
for some $K > 100$ to be determined. We take $K$ much larger than the $K_1$ from Lemma \ref{lem:av-ll-1} so that $\tau_2 > S_1$ with probability at least $1- N^{-D}$. 

Compared to the proof of Lemma \ref{lem:av-ll-1} we  first have that $N\eta_t\rho_t\ge (\log N)^3$ for $0\le t\le S_2$,  
instead of $N\eta_t\rho_t\ge N^{5\xi}$. This will be compensated with the fact that various $N^{\xi}$ factors in the estimates are replaced by logarithmic factors due to the definition of $\tau_2$. 

The estimates \eqref{eq:flowgro} and \eqref{eq:flowgro-2} become
\begin{align}
\label{eq:flowgro-3}
\big|X_{1}^\av (t \wedge \tau_t ) \big|&\leq \big|p_{S_1,t\wedge\tau_2}X_{1}^\av (S_1) \big|+\left|\int_{S_1}^{t\wedge\tau_2}p_{s,t\wedge\tau_1}\,\d N_s\right|+\int_{S_1}^{t\wedge\tau_2} p_{s,t\wedge\tau_2}\big|X_{1}^\av (s) X_{2}^\av (s) \big|\,\d s \notag \\
&\quad+4\int_{S_1}^{t\wedge\tau_1} p_{s,t\wedge\tau_2}\frac{\rho_s}{N\eta_s^2}\,\d s, 
\end{align}
and
\begin{align}\label{eq:flowgro-4}
\big|X_{2}^\av (t \wedge \tau_2 ) \big|&\leq \big|p_{S_1,t\wedge\tau_2}^2X_{2}^\av (S_1) \big|+\left|\int_{S_1}^{t\wedge\tau_1}p_{s,t\wedge\tau_2}^2\,\d \widehat{N}_s\right|+\int_{S_1}^{t\wedge\tau_1} p_{s,t\wedge\tau_1}^2\big| X_{2}^\av (s) \big|^2\,\d s \notag \\
&\quad+2\int_{S_1}^{t\wedge\tau_2} p_{s,t\wedge\tau_2}^2\left(\frac{\rho_s}{\eta_s^2}\big|X_{1}^\av (s) \big|+4\frac{\rho_s}{N\eta_s^3}\right)\,\d s.
\end{align}
First, by applying Lemma \ref{lem:stochastic} the analog of \eqref{eqn:aaa-2} holds, with the suprema being over $u \in [S_1, S_2]$ (instead of $u \in [0, S_1]$), and taking the $K$ here in the definition of $\tau_2$ in \eqref{eqn:tau-2-def} to be large enough so that the analog of \eqref{eqn:aaa-2} holds with this choice of $K$. 

In a similar fashion to \eqref{eq:usbound} we have, 
\beq
\begin{split}
\int_{S_1}^{t\wedge\tau_2} p_{s,t\wedge\tau_2}\big|X_{1}^\av (s) X_{2}^\av  (s) \big|\,\d s&\leq C K^2\int_{S_1}^{t\wedge\tau_2} \frac{\eta_s}{\eta_{t \wedge\tau_2}}\frac{\log N}{N\eta_s}\frac{( \log N)^2}{N\eta_s^2}\leq \frac{C K^2}{N\eta_{t \wedge \tau_2}}, \\
\int_0^{t\wedge\tau_1}p_{s,t\wedge\tau_1}\frac{\rho_s}{N\eta_s^2}\,\d s&\le 2 \int_0^{t\wedge\tau_1}\frac{\rho_s}{N\eta_s\eta_{t\wedge\tau_1}}\,\d s\leq 3 \frac{\log N}{N\eta_{t\wedge\tau_1}}.
\end{split}
\eeq
Above we used $N \rho_s \eta_s \geq ( \log N)^3$ for $s \leq S_2$.  
Similar to \eqref{eq:usbound2} we have,
\beq
\begin{split}
\int_{S_1}^{t\wedge\tau_2} p_{s,t\wedge\tau_1}^2\big| X_{2}^\av (s) \big|^2\,\d s&\leq C K^2 \int_{S_1}^{t\wedge\tau_2}\frac{\eta_s^2}{\eta_{t\wedge\tau_2}^2}\frac{( \log N)^4}{N^2\eta_s^4}\,\d s\leq CK^2 \frac{\log N}{N\eta_t^2}, \\
\int_{S_1}^{t\wedge\tau_2} p_{s,t\wedge\tau_2}^2\left(\frac{\rho_s}{\eta_s^2}\big|X_{1}^\av (s) \big|+4\frac{\rho_s}{N\eta_s^3}\right)\,\d s &\leq \int_{S_2}^{t\wedge\tau_2}\frac{\eta_{s}^2}{\eta_{ t\wedge\tau_2}^2}\left(\frac{6 K \rho_s\log N}{N\eta_s^3}\right)\,\d s\\
&\le 7 K \frac{(\log N)^2}{N\eta_{t\wedge\tau_2}^2}.
\end{split}
\eeq


The initial data terms are estimated using \eqref{eq:propag} and \eqref{eq:imp-G-1} by,
\beq
\begin{split}
|p_{S_1,t\wedge\tau_2}X_{1}^{\av} (S_1)|&\le 2\frac{\eta_{S_1}}{\eta_{t \wedge \tau_2}} \frac{K_1 \log N}{N\eta_{S_1}}\le \frac{2 K_1 \log N}{N\eta_{t\wedge\tau_2}}, \\
|p_{S_1,t\wedge\tau_2}^2X_{2}^{\av} (S_1) |&\le 4\frac{\eta_{S_1}^2}{\eta_{t \wedge \tau_2}^2} \frac{K_1(\log N)^2}{N\eta_{S_1}^2}\le \frac{4 K_1 (\log N)^2}{N\eta_{t\wedge\tau_2}^2} .
\end{split}
\eeq
with probability at least $1- N^{-D}$. 
Therefore, in an analogous fashion to the proof of Lemma \ref{lem:av-ll-1} 
we conclude, after taking $K$ much larger than $K_1$, that 
\beq
|X_{1}^{\av} (t \wedge \tau_2) |\le \frac{4 K\log N}{N\eta_{t\wedge\tau_2}},
\qquad\quad  |X_{2}^{\av} (t\wedge\tau_2)|\le \frac{19 K (\log N)^2}{N\eta_{t\wedge\tau_2}^2}.
\eeq
This implies $\tau_2=S_2$ with probability at least $1 - 3 N^{-D}$ and completes the proof of the lemma. \qed 

\bel \label{lem:av-ll-3}
For any $D>0$ there is a $K_3 = K_3 (D)>0$ so that with probability at least $1- N^{-D}$ we have for all $t \in [S_2, T]$ that
\beq \label{eqn:av-ll-3}
|X_1^\av (t) | \leq \frac{ K_3 \log N}{N \eta_t}. 
\eeq
\eel
\proof Unlike in Lemmas \ref{lem:av-ll-1} and \ref{lem:av-ll-2} we only consider the evolution of $\langle G_t\rangle$.  Recall the short-hand notation \eqref{eq:shortnot} and define the stopping time
\beq
\tau_3:=\inf \left\{t\ge S_2 : |X_{1}^{\av} (t) |=20 \frac{K \log N}{N\eta_t}\right\}\wedge T.
\eeq
We will take $K$ much larger than the $K_1, K_2$ from Lemmas \ref{lem:av-ll-1} and \ref{lem:av-ll-2}. 
We point out that for $0\le t\le T$ we have $N\eta_t\rho_t\ge K_0 \log N$, for $K_0=100 K^2$. 

We now recall the equation describing the evolution of the single resolvent $G_t$:
\beq
\label{eq:neweqsg}
\d \langle G_t-M_t\rangle =\d N_t+\langle G_t-M_t\rangle\left(\frac{1}{2}+\langle G_t^2\rangle\right) \d t+ \frac{\1_{ \{ \beta=1\}}}{N} \tilde{\sum}_{ij} \langle G_t^2 E_i G_t^\mft E_j \rangle \d t.
\eeq
While \eqref{eq:neweqsg} contains $\langle G_t^2\rangle$, in the following we will see that we can now afford to use the estimate
\beq
\label{eq:newsw}
\big|\langle G_t^2\rangle\big|\le \frac{\rho_t}{\eta_t}+20 \frac{K \log N}{N\eta_t^2} =: r_s
\eeq
for $t\le\tau_3$, and so that it is not necessary to monitor the evolution of $\langle G_t^2\rangle$. This is a consequence of the fact that in this last part we change the scale on which the local law holds just by a power of $\log N$.

We start with the estimate of the martingale in \eqref{eq:neweqsg}. The quadratic variation process is bounded by \eqref{eq:boundqv}, and so by
 \eqref{eqn:sub-gaussian} and an easy discretization argument (see (3.36)-(3.37) of \cite{cipolloni2024maximum}) we have,
\beq \label{eqn:ll-3-stochastic}
\pp\left[\sup_{S_2 <u < \tau_3 } N\eta_u\big|N_{u}\big|> K \sqrt{\log N}\right]\le N^{-D}
\eeq
after taking $K$ large enough depending on $D$. 
Integrating \eqref{eq:neweqsg} and using \eqref{eq:newsw}, that $|\langle G_t^2E_iG_t^\mathfrak{t}E_j\rangle|\le 2\rho_t/\eta_t^2$ by \eqref{eqn:schwarz-ave}, and \eqref{eqn:ll-3-stochastic} we obtain
\beq
|X_{1}^{\av} (t \wedge \tau_3) |\le |X_{1}^{\av}(S_2) |+\int_{S_2}^{t\wedge\tau_3} r_s |X_{1}^\av (s)|\,\d s+2\frac{K\sqrt{\log N}}{N\eta_{t\wedge\tau_3}}.
\eeq
Since $\rho_T \eta_T N \geq 100 K^2  \log N$ we have the inequalities,
\beq
\exp \left( \int_{s}^t r_u \d u  \right) \leq 2 \frac{\eta_s }{\eta_t} , \qquad r_u \leq 2 \frac{ \rho_u}{\eta_u}
\eeq
as long as $K$ is large enough. Gr\"onwall's inequality in integral form gives,
\begin{align}
X_1^{\av} (t \wedge \tau_3) \leq & |X_{1}^{\av}(S_2) | + 2\frac{K\sqrt{\log N}}{N\eta_{t\wedge\tau_3}} \notag\\
+ & \int_{S_2}^{ t \wedge \tau_3} r_s \left( |X_{1}^{\av}(S_2) | + 2\frac{K\sqrt{\log N}}{N\eta_{s\wedge\tau_3}}  \right)  \exp \left( \int_{s}^t r_u \d u  \right)  \d s .
\end{align}

Therefore with probability at least $1 - 2 N^{-D}$, using Lemma \ref{lem:av-ll-2} we have, by taking $K$ much larger than $K_2$,
\beq
\begin{split}
|X_{1,t\wedge\tau_3}^{\av}|&\le \frac{K_2 \log N}{N\eta_{S_2}}+ 2\frac{K\sqrt{\log N}}{N\eta_{t\wedge\tau_3}}+4 \int_{S_2}^{t\wedge\tau_3}\frac{\eta_s}{\eta_{t\wedge\tau_3}}\cdot\frac{\rho_s}{\eta_s}\cdot \left(\frac{K_2\log N}{N\eta_{S_2}}+2\frac{K\sqrt{\log N}}{N\eta_s}\right)\, \d s \\
&\le \frac{K \log N}{N\eta_{S_2}}+\frac{9 K \sqrt{\log N}\log\log N}{N\eta_{t\wedge\tau_3}}+\frac{4K_2\rho_T(T-S_2)\log N}{N\eta_{S_2}\eta_{t\wedge\tau_3}} \\
&\le 4 K \frac{\log N}{N\eta_{t\wedge\tau_3}}.
\end{split}
\eeq
We point out that in the second inequality we used 
\beq
\log ( \eta_{S_2} / \eta_T ) \leq C \log \log N
\eeq
by Lemma \ref{lem:chars}, 
and that in the last inequality we used $t\wedge\tau_3-S_2\le T-S_2\le 2\eta_{S_2}/\rho_{S_2}$ (the last inequality following from \eqref{eqn:del-t-ratio})and $\rho_T\le \rho_{S_2}$. This shows that $\tau_3=T$ with probability at least $1-2 N^{-D}$ and thus concludes the proof of the Lemma. \qed 

\vspace{5 pt}

\noindent{\bf Proof of Proposition \ref{prop:ll-id}}. This follows immediately from Lemmas \ref{lem:av-ll-2} and \ref{lem:av-ll-3}. \qed

\smallskip

\bel \label{lem:B-flow}
Let $D>0$. There is a $K = K(D)$ so that with probability at least $1- N^{-D}$ it holds for all $t \in [S_1, T]$ that,
\beq
\left| \langle (G_t - M_t ) B \rangle \right| \leq K \frac{ \| B \| \log N}{N \eta_t}
\eeq
\eel
\proof To simplify the notation we assume $\| B \| =1 $. 
With Proposition \ref{prop:ll-id} in hand, the averaged local law for general observables is substantially simpler. In fact, we will now see that it is not needed to study the evolution of $G_t^2$ and no Gr\"onwall inequality and stopping times are needed.

Similarly to \eqref{eq:flownocharchar}, we obtain
\beq
\label{eq:neweqsgB}
\d \langle (G_t-M_t)B\rangle =\d N^B_t+\frac{1}{2}\langle (G_t-M_t)B\rangle \d t+\langle G_t^2B\rangle\langle G_t-M_t\rangle \d t+ \frac{\1_{ \{ \beta=1\}}}{N} \tilde{\sum}_{ij} \langle G_t^2B E_i G_t^\mft E_j \rangle \d t,
\eeq
where
\beq
\d N^B_t:=-\frac{1}{\sqrt{N}}\langle G_tBG_t\d \mathfrak{B}_t\rangle.
\eeq
Here $\mathfrak{B}_t$ is defined in \eqref{eq:defbigBM}. Proceeding similarly to the case $B=I$ it is easy to see that
\beq
\d [ \widebar{N}^B_t , N^B_t]=\frac{1}{N^2\eta_t^2}\langle\Im G_tB\Im G_tB^*\rangle\d t+\frac{\1_{ \{ \beta=1\}}}{N^2}\langle (G_tBG_t)(G_t^*BG_t^*)^\mathfrak{t}\rangle\d t \le 4  \frac{\rho_t}{N^2 \eta_t^3} \d t
\eeq
(with the inequality being proved similarly to Lemma \ref{lem:schwarz-ave}) 
and so, in the same manner as in the proof of Lemma \ref{lem:av-ll-3} we have
\beq \label{eqn:B-stoch}
\pp\left[ \sup_{0 < u < T} (N \eta_u) |N_u^B| > K \sqrt{ \log N} \right] \leq N^{-D}. 
\eeq
for some $K=K(D) \geq 100$ depending on $D$. 
Integrating \eqref{eq:neweqsgB} and using $N^{-1} | \langle G_t^2B E_i G_t^\mft E_j \rangle | \leq 2 \frac{ \rho_t}{N \eta_t^2}$ (from \eqref{eqn:schwarz-ave}) and \eqref{eqn:B-stoch} we see that,
\beq \label{eqn:ll-B-1}
\big|\langle (G_t-M_t)B\rangle\big|\le e^{t/2}\big|\langle (G_0-M_0)B\rangle\big|+\left| \int_0^t \e^{(t-s)/2}\langle G_s^2B\rangle\langle G_s-M_s\rangle \d s \right| +2\frac{K \sqrt{\log N}}{N\eta_t} .
\eeq
From \cite[Theorem 5.2]{alt2021spectral} and \eqref{eqn:chars-1} we have with overwhelming probability,
\beq \label{eqn:inB}
\big|\langle (G_0-M_0)B\rangle\big| \leq \frac{ N^{\xi}}{N \eta_0} \leq \frac{1}{N \eta_t}. 
\eeq
for any $t\in [S_1, T]$. It remains to estimate the integral in \eqref{eqn:ll-B-1}. Using $|\langle G_t^2 B\rangle|\le 2 (\rho_t\lVert B\rVert)/\eta_t$, \eqref{eq:impsingG}, and \eqref{eqn:prop-ll-av} we have for any $t \in [S_1, T]$ that (after possibly increasing the value of $K$)
\begin{align} \label{eqn:B-ll-a1}
    & \int_0^t |\langle G_s^2 B \rangle \langle G_s - M_s \rangle | \d s \leq K \frac{ \eta_0}{\kappa_0} \int_0^{S_1} \frac{ N^{\xi} \rho_s}{N \eta_s \eta_{s \wedge t_*}} + K \int_0^t \frac{ \rho_s \log N}{N \eta_s^2} \d s  \notag\\
    \leq & K \frac{\eta_0}{ \kappa_0 \eta_{S_1 \wedge t_*} } \int_0^{S_1} \frac{ N^{\xi} \rho_s}{N \eta_s} \d s + K \frac{ \log N}{N \eta_s} \leq 2 K \frac{ \log N}{N \eta_t} 
\end{align}
with probability at least $1- 2 N^{-D}$. In the last inequality we used \eqref{eq:ei-1}. This completes the proof. \qed



\bel \label{lem:iso-flow}
For any $D>0$ there is a $K>0$ so that for any two unit vectors ${\bf x}$ and ${\bf y}$ we have with probability at least $1- N^{-D}$ that
\beq
|(G_t -M_t)_{{\bf x}, {\bf y}}| \leq K \sqrt{ \frac{ \rho_t \log N}{N \eta_t}}
\eeq
for all $t \in [S_1, T]$.
\eel
\proof 
For unit vectors ${\bf u}, {\bf v}\in\cc^{2N}$, we introduce the short hand-notation
\beq
X_{{\bf u}, {\bf v}}^{\iso} (t):=(G_t-M_t)_{{\bf u}{\bf v}}.
\eeq 
Next, we define the stopping time,
\beq
\tau:=\inf \left\{t\geq 0 : \sup_{{\bf v},{\bf w}\in \{{\bf x},{\bf y}\}}|X_{{\bf v},{\bf w}}^{\iso} (t)|>\rho_t\right\}\wedge T.
\eeq
In the remainder of the proof the vectors ${\bf u}, {\bf v}$ are arbitrary vectors in $\{ {\bf x}, {\bf y}\}$. 
For $t \in [0, S_1]$ we have $N \eta_t \rho_t \geq N^{5 \xi}$ and so by \cite[Theorem 5.2]{alt2021spectral}, we have
\beq
\label{eq:iniso}
|X_{{\bf u}, {\bf v}}^{\iso}(t) |\leq N^\xi \sqrt{\frac{\rho_t}{N\eta_t}} \leq N^{-\xi} \rho_t,
\eeq
with overwhelming probability, which implies $\tau >S_1$. To estimate $X_{{\bf u}, {\bf v}}^{\iso}(t)$ for  $ t \in [S_1, T]$ we integrate \eqref{eq:flownocharisochar} in time, using \eqref{eqn:schwarz-iso} 
we obtain
\beq \label{eqn:iso-proof-a1}
\big|X_{{\bf u}, {\bf v}}^\iso (t \wedge \tau) \big|\lesssim \big|X_{{\bf u}, {\bf v}} (0)\big|+\big|\widetilde{N}_{t\wedge\tau}\big|+2\int_0^{t\wedge\tau}\left(\big|X_{1}^\av (s)\big|\frac{\rho_s}{\eta_s}+2\frac{\rho_s}{N\eta_s^2}\right)\,\d s 
\eeq
The second term in the integrand integrates to $\frac{2}{N \eta_{t \wedge \tau}}$. The first term in the integral is bounded in the exact same fashion as \eqref{eqn:B-ll-a1} and so,
\beq \label{eqn:iso-proof-a2}
2\int_0^{t\wedge\tau}\left(\big|X_{1}^\av (s)\big|\frac{\rho_s}{\eta_s}+2\frac{\rho_s}{N\eta_s^2}\right)\,\d s  \leq K \frac{ \log N}{N \eta_{t \wedge \tau}}
\eeq
with probability at least $1- N^{-D}$ by possibly increasing the value of $K$.


We now estimate the martingale term. We have for $s < \tau$,  
\begin{equation}
\label{eq:boundqviso}
\d [\overline{\widetilde{N}_s},\widetilde{N}_s]=\frac{1}{N\eta_s^2}(\Im G_s)_{\bf u u}(\Im G_s)_{\bf v v}+\frac{ \1_{ \{ \beta=1\}}}{N}(G_s(G_s^*)^\mft)_{\bf u u}(G_s(G_s^*)^\mft)_{\bf v v} \le 4 \frac{\rho_s^2}{N\eta_s^2},
\end{equation}
with the inequality following similarly to the proof of Lemma \ref{lem:schwarz-iso} 
and so, similarly to \eqref{eqn:aaa-2}, that
\beq
\label{eqn:aaa-2iso}
\pp\left[ \sup_{ 0 < u<s } \sqrt{\frac{N\eta_{u\wedge\tau}}{\rho_{u\wedge\tau}}} | \widetilde{N}_{u\wedge \tau}| >  K \sqrt{\log N} \right]\le N^{-D},
\eeq
for  $K$ sufficiently large depending on $D$. 

By combining \eqref{eq:iniso}, \eqref{eqn:iso-proof-a1}, \eqref{eqn:iso-proof-a2}, and \eqref{eqn:aaa-2iso}, (as well as \eqref{eqn:chars-1}) we obtain for $t \in [S_1, T]$,
\beq
|X_{{\bf u}, {\bf v}}^{\iso} (t \wedge \tau) |\le C K\sqrt{\log N} \sqrt{\frac{\rho_{t\wedge\tau}}{N\eta_{t\wedge\tau}}}. 
\eeq
By the assumption that $\rho_T \eta_T N \geq K_0 \log N$, this shows us that $\tau = T$, as the RHS can be made to be less than $\frac{\rho_{t\wedge \tau}}{2}$ after increasing the value of $K_0$ depending on $K$. Then, the above estimate is what we wanted to prove. \qed 

\appendix

\section{Conclusion of the Proof of Theorem \ref{theo:isolaw}}
\label{sec:conclproof}

In this section we pick up where we left off in the proof of Theorem \ref{theo:isolaw}.  Let us summarize what we have shown so far: for any $D>0$ there is a constant $K = K (D)$ (depending only on the constants in Definition \ref{def:model} and $D$) so that the following holds. If $\xi >0$ and $Y$ is a matrix of the form $Y = (1-T)^{1/2} W + \sqrt{T} Z$, with $Z$ a matrix of i.i.d. Gaussians and $W$ as in Definition \ref{def:model}, and $T := N^{-\xi}$, then for any $z$ and $\eta$ satisfying $ K N^{-1} \log N \leq \eta \rho^z ( \i \eta) \leq N^{ - 100 \xi} $ and $|z| \leq 1 + N^{-100 \xi}$ we have that for any fixed $B$ and unit vectors ${\bf x}, {\bf y}$ that
\beq
\langle (G_Y^z ( \i \eta ) - M^z ( \i \eta ))B \rangle \leq \frac{K \| B \| \log N}{N \eta},
\eeq
and
\beq
\langle {\bf x}, (G_Y^z ( \i \eta ) - M^z ( \i \eta ) ) {\bf y} \rangle \leq K \sqrt{ \log N } \sqrt{ \frac {\rho^z ( \i \eta)}{N \eta}},
\eeq
with probability at least $1 - N^{-D}$. Here we denoted $G_Y^z ( \i \eta)$ to be the resolvent of the Hermitization of $Y-z$ as defined in \eqref{eq:defres}, but with $X= Y$. The above holds for $N$ large enough depending on $\xi >0$ and the constants in Definition \ref{def:model}. 

We now let $X$ be a general i.i.d. matrix as in Definition \ref{def:model}. We wish to show that the above estimates hold for $Y$ replaced by $X$. By a standard construction (see \cite[Lemma 16.2]{erdHos2017dynamical})  we can choose $W$ so that with $Y = (1- T)^{1/2} W + \sqrt{T} Z$ as above we have,
\beq
\left| \ee[X_{ij}^a \bar{X}_{ij}^b] - \ee[Y_{ij}^a \bar{Y}_{ij}^b] \right| \leq 0+ N^{-2}\delta_{\{ a+b=4\}}  T, \qquad 0 \leq a + b \leq 4.
\eeq
Let now $Z_i(A)$ denote the functions on the space of $N \times N$ matrices given by,
\beq
Z_1(A) = N \eta \langle (G^z_A ( \i \eta) - M^z ( \i \eta ) )B \rangle \qquad  Z_2(A) = (N \eta / \rho^z ( \i \eta ) )^{1/2} \langle {\bf x} , (G^z_A ( \i \eta) - M^z ( \i \eta ) ) {\bf y} \rangle .
\eeq
Let $\Phi$ be an enumeration of the matrix indices $[ N ]^2$ of $N \times N$ matrices (i.e., $\Phi : [  N ]^2 \to [ N^2 ]$  is a bijection). For $j \in [ N^2]$ let $Y^{(j)}$ be the matrix whose entries are $Y^{(j)}_{ab} = X_{ab}$ if $\Phi(a, b) \leq j$ and $Y^{(j)}_{ab} = Y_{ab}$ if $\Phi (a, b) > j$. Note that $Y^{(0)} = Y$ and $Y^{(N^2)} = X$. 
Let us now define,
\beq \label{eqn:pk-def}
p_k^{(n)} := \sup_{0 \leq j \leq N^2 } \pp\left[ |Z_n(Y^{(j)}) | > K (\log N)^{1/n} + k \right].
\eeq
In the remainder of the proof we will show that for any $L>0$ there is a $C_2 >0$ so that for all $k \geq 2$ we have that
\beq \label{eqn:a-comp}
p_k^{(n)} \leq \pp\left[ |Z_n (Y)| > K ( \log N)^{1/n} \right] + C_2 p_{k-1}^{(n)} T^{1/2} + C_2 N^{-L}. 
\eeq
By iterating the above estimate an $N$-independent number of times (depending on $D/ \xi$) one obtains the estimates of  Theorem \ref{theo:isolaw} for the matrix $X$, at the cost of increasing the constant $K$.

It remains to give the proof of \eqref{eqn:a-comp}. For notational simplicity we will consider only the case of real symmetric matrices, the complex case being similar. The proof of this estimate is largely similar to the proof of \cite[Proposition 3.7]{cipolloni2024maximum} so we will not provide all of the details. Let us first consider the case $n=1$ and $k \geq 2$. Let $F : \cc \to \rr$ be a smooth function so that $F(z) =1$ for $|z| > C_1 \log N + k -1/4$ and $F(z) = 0$ for $|z| < C_1 \log N + k -3/4$. We may assume that $\| F \|_{C^5} \leq C$ uniformly in $k$. 

We need to estimate one of the terms on the right-hand side of \eqref{eqn:pk-def} for fixed $j$. In order to do so we write,
\begin{align} \label{eqn:gfct-a3}
& \pp\left[ |Z_1 (Y^{(j)} ) | > C_1 \log N + k \right]  \leq \ee[ F(Z_1(Y^{(j)} )] \notag\\
&\qquad\qquad\quad=  \ee[ F (Z_1( Y^{(0)} ) )] + \sum_{m=1}^j \ee[ F (Z_1 ( Y^{(m)} ) )] - \ee[ F (Z_1 ( Y^{(m-1)} ) )] .
\end{align}
The first term on the second line is bounded by the first term on the RHS of \eqref{eqn:a-comp}. Fix now a term $m$ in the sum on the RHS and assume $\Phi (a,b) = m$. Let us introduce the notation,
\beq
\hata := a + N, \quad \hatb = b + N . 
\eeq
For $\theta \in \rr$, let $W(\theta)$ be the matrix such that $W(\theta)_{ab} = \theta$ and $W(\theta)_{ij} = Y^{(m)}_{ij} = Y^{(m-1)}_{ij}$ otherwise. By Taylor expansion to fifth order we have for any $\eps >0$,
\begin{align} \label{eqn:gfct-a1}
& \left|  \ee[ F (Z_1 ( Y^{(m)} ) )] - \ee[ F (Z_1 ( Y^{(m-1)} ) )]  \right| \notag \\
&\leq  N^{-L-10} +N^{\eps}(TN^{-2}+N^{-5/2} ) \ee\left[ \sup_{| \theta'| \leq N^{\eps-1/2}} \sup_{m\leq 5} \left| \del_\theta^m  F(Z_1 (W(\theta) )) \vert_{\theta=\theta'} \right| \right]
\end{align}
In preparation for estimation of the derivatives on the RHS we first note that by resolvent expansion one can obtain the estimates, 
\beq \label{eqn:gfct-a2}
\sup_{ |\theta| \leq N^{\eps-1/2}} | \langle \be_i , G^z_{W(\theta)} ( \i \eta)  \be_j \rangle -\langle \be_i , M^z ( \i \eta) \be_j \rangle | \leq N^{2\eps} \sqrt{ \frac{ \rho^z(\i \eta)}{N \eta}}
\eeq
and
\begin{align} \label{eqn:a-iso-2}
 (G^z_{W(\theta)} )_{ij} - (G_{Y^{(m)}}^z)_{ij} &= ( \theta - Y_{ab} ) \left( (G^z_{Y^{m}} )_{ia} (G^z_{Y^{m}} )_{\hatb j} +  (G^z_{Y^{m}} )_{i \hatb} (G^z_{Y^{m}} )_{a j} \right) \notag\\
+ & (\theta - Y_{ab} )^2 \bigg\{  (G^z_{Y^{m}} )_{i a}  (G^z_{Y^{m}} )_{\hatb a}  (G^z_{Y^{m}} )_{\hatb j} + (G^z_{Y^{m}} )_{i a}  (G^z_{Y^{m}} )_{\hatb \hatb}  (G^z_{Y^{m}} )_{a j}  \notag\\
 & \quad + (G^z_{Y^{m}} )_{i \hatb}  (G^z_{Y^{m}} )_{a a}  (G^z_{Y^{m}} )_{\hatb j} + (G^z_{Y^{m}} )_{i \hatb}  (G^z_{Y^{m}} )_{a \hatb}  (G^z_{Y^{m}} )_{a j} \bigg\} \notag\\
+ &\O ( N^{5 \eps-3/2} ) \sum_{u,v \in \{a, \hatb \}}  |(G^z_{Y^{m}} )_{iu} (G^z_{Y^{m}} )_{vj}| + \O (N^{-5})
\end{align}
with overwhelming probability. Note that the second estimate, together with the fact that $\| B \| \leq C$, imply that for $n=1, 2$, with overwhelming probability,
\beq
\sup_{ |\theta| \leq N^{\eps-1/2}} |Z_n (W(\theta) ) - Z_n (Y^{(m-1)} ) | \leq N^{-\xi/10},
\eeq
by taking $\eps >0$ small enough depending on $\xi$. We also used that $\eta / \rho \leq N^{-\xi}$ by assumption  in the case $n=2$. On the event that the above estimate holds, we see that the derivative in the expectation in \eqref{eqn:gfct-a1} is non-zero only if $|Z_1 (Y^{(m-1} ) | > C_1 \log N + k-1$, by the support properties of $F$. Hence, in order to deduce \eqref{eqn:a-comp} in the case $n=1$ we see that it suffices to show that with overwhelming probability,
\beq \label{eqn:Z1-bd}
\sup_{| \theta'| \leq N^{\eps-1/2} } \left| \del_\theta^m  Z_1 (W(\theta) )) \vert_{\theta=\theta'} \right| \leq N^{C \eps} ,
\eeq
for some $C>0$. 
By direct calculation, the derivative on the left-hand side is a finite sum of terms of the form
\beq
N \eta \frac{1}{N} \sum_{ij} B_{ij} (G^z_{W(\theta)} )_{il} (G^z_{W(\theta)} )_{kj} P,
\eeq
where: $P$ is a monomial in resolvent entries not depending on $i$ or $j$; $l$ and $k$ are fixed indices that are either $a$, $b$ or $a+N$ or $b+N$. The desired estimate then follows by the Ward identity, \eqref{eqn:gfct-a2} and that $\| B \| \leq C$. This completes the proof of \eqref{eqn:a-comp} when $n=1$. The case $n=2$ is similar. The only change (apart from changing $\log N$ to $ \sqrt{ \log N}$ in the choice of support of $F$) is that \eqref{eqn:Z1-bd} is replaced by,
\beq \label{eqn:a-iso-3}
\sup_{| \theta'| \leq N^{\eps-1/2} } \left| \del_\theta^m  Z_1 (W(\theta) )) \vert_{\theta=\theta'} \right| \leq N^{5\eps} \sqrt{ \frac{ N \eta}{ \rho}} \left( \sum_{ l \in S_{ab} } x_l^2 + y_l^2 + \frac{ \rho}{N \eta} \right)
\eeq
where $S_{ab} = \{ a, b, a+N, b+N\}$  and $x_l,y_l$ denote the entries of the unit vectors ${\bf x}, {\bf y}\in\cc^{2N}$. The proof of the above estimate follows below. We note that \eqref{eqn:a-iso-3} is sufficient to complete the proof, because its overall contribution to the telescoping sum in \eqref{eqn:gfct-a3} is,
\begin{align}
& p_{k-1}^{(2)} \sum_{a,b} \left\{(T N^{-2} + N^{-5/2} ) \left( \sum_{ l \in S_{ab} } x_l^2 + y_l^2 + \frac{ \rho}{N \eta} \right) \right\} \sqrt{\frac{N \eta}{\rho}} \notag\\
\leq & p_{k-1}^{(2)} N^{C \eps} T\left(  \sqrt{ \frac{\rho}{N \eta}} + \sqrt{ \frac{\eta}{\rho N}}  \right) \leq p_{k-1}^{(2)} N^{C \eps} T,
\end{align}
using that ${\bf x}$ and ${\bf y}$ are unit vectors. 

In order to prove \eqref{eqn:a-iso-3}, we note that the derivatives are finite sums of terms of the form, 
\beq \label{eqn:a-iso-1}
\sqrt{ \frac{N\eta}{\rho}} P \sum_{ij} x_i (G^z_{W(\theta)} )_{il} (G^z_{W(\theta)} )_{mj} y_j
\eeq
where $m, l \in S_{ab}$ and $P$ is a  monomial of resolvent entries. We need to show that \eqref{eqn:a-iso-1} is bounded by the RHS of \eqref{eqn:a-iso-3}.  

First note that the expansion \eqref{eqn:a-iso-2} implies that for any unit vector,
\beq
\sup_{ |\theta| \leq N^{\eps-1/2}} \left| \langle {\bf v} , (G^z_{W(\theta)} - G_{Y^{(m)}} ) {\bf e}_k \rangle \right| \leq \frac{ N^{5\eps}}{\sqrt{N}}
\eeq
with overwhelming probability. By \cite[Theorem 3.1]{cipolloni2021fluctuation}, \cite[Theorem 3.1]{cipolloni2024precise} used in the regimes $1-|z|\ge c$ and $|1-|z||<c$, respectively, we have 
\beq
|  \langle {\bf v}, G_{Y^{(m)}} {\bf e}_k \rangle - \langle {\bf x} , M {\bf e}_k \rangle | \leq N^{\eps} \sqrt{ \frac{ \rho}{N \eta}}
\eeq
with overwhelming probability
From the above two estimates,
\beq
\sup_{ |\theta| \leq N^{\eps-1/2}} | \langle {\bf v}, G^z_{W(\theta)} {\bf e}_k \rangle - \langle {\bf v} , M {\bf e}_k \rangle | \leq N^{4 \eps} \sqrt{ \frac{ \rho}{N\eta}}
\eeq
for ${\bf v} = {\bf x}, {\bf y}$. 
By applying this estimate to the terms \eqref{eqn:a-iso-1} one proves the estimate \eqref{eqn:a-iso-3}. \qed 

\section{Stochastic terms}

\bel \label{lem:stochastic}
Let $T\leq 1$ and let $( \eta_t, z_t)_{t \in[0, T]}$ be a characteristic satisfying $|z_t| \leq 10$ and $\eta_t \geq N^{-1}$ for all $t$. Let $\tau$ be a stopping time so that for $t < \tau$ we have that
\beq
\langle \Im[ G_{t}^{z_t} ( \i \eta_t)] \rangle \leq 10 \rho_t
\eeq
almost surely. For any $D>0$ there is a $C(D) >0$ so that the following holds. For any $0 \leq T_1 < T_2 \leq T$ we have that,
\beq
\pp\left[ \exists t \in [T_1, T_2] : \left| \int_{T_1 \wedge \tau}^{t \wedge \tau} p_{s, t} \d N_s \right| > C \frac{ \left( \log ( \eta_{T_1} / \eta_{T_2} ) \log N \right)^{1/2}}{N \eta_t} \right] \leq N^{-D}
\eeq
and
\beq
\pp\left[ \exists t \in [T_1, T_2] : \left| \int_{T_1 \wedge \tau}^{t \wedge \tau} p^2_{s, t} \d \hat{N}_s \right| > C \frac{ \left( \log ( \eta_{T_1} / \eta_{T_2} ) \log N \right)^{1/2}}{N \eta_t^2} \right] \leq N^{-D}
\eeq
\eel
Before embarking on the proof we first recall the standard fact that if $M_t$ is a Martingale whose quadratic variation is bounded almost surely by some deterministic constant $B$ then
\beq \label{eqn:sub-gaussian}
\pp\left[ \sup_{t} |M_t | > 100 s B^{1/2} \right] \leq 2 \e^{ - s^2}. 
\eeq
Additionally, we have for the quadratic variation of $N_s$ that, on the event that $| \langle G_t -M_t \rangle| \leq \rho_t$ that,
\begin{align}
\label{eq:boundqv}
\d [N_s,\bar{N}_s] = & \tilde{\sum}_{ij} \left( \frac{1}{N^2}\langle G_s(w_s)^2E_iG_s(\overline{w_s})^2E_j\rangle + \frac{\1_{ \{ \beta=1\}} }{N^2} \langle G_s(w_s)^2 E_i [ G_s (\bar{w}_s)^2]^\mft E_j \rangle \right) \d s \notag\\
\le & \frac{8\rho_s}{N^2\eta_s^3}\d s.
\end{align}
with the inequality following in a similar manner to the proof of Lemma \ref{lem:schwarz-ave}. By a similar calculation, 
\beq
\label{eq:boundqv2}
\d [\widehat{N}_s,\overline{\widehat{N}}_s]\le \frac{20\rho_s}{N^2\eta_s^5}\,\d s.
\eeq
\noindent{\bf Proof of Lemma \ref{lem:stochastic}}.  We only prove the first estimate, the second being similar. For integer $0 \leq k \leq N^{10}$ define $t_k:= T_1 + k N^{-10} (T_2-T_1)$. It is easy to see that for $t \in [t_k, t_{k+1}]$ that,
\beq \label{eqn:a-eta}
\frac{1}{2} \eta_{t_k} \leq \eta \leq 2 \eta_{t_k} .
\eeq
Moreover, for $t \in [t_k, t_{k+1}]$ we have, since $p_{s,t} = p_{s, t_k} p_{t_k, t_{k+1}}$
\beq
\left| \int_{T_1 \wedge \tau} ^t p_{s, t} \d N_s \right| \leq 2 \left| \int_{T_1 \wedge \tau} ^t p_{s, t_k} \d N_s \right|
\eeq
almost surely. By the assumption on $\tau$ we have for $t < \tau$ by \eqref{eq:boundqv}
\beq
\left\langle \int_{T_1 \wedge \tau} ^t p_{s, t_k} \d N_s  , \int_{T_1 \wedge \tau} ^t p_{s, t_k} \d \bar{N}_s \right\rangle \leq \frac{8}{(N \eta_{t_k})^2} \log ( \eta_{T_1} / \eta_{T_2} )
\eeq
Therefore, by \eqref{eqn:sub-gaussian} we have for sufficiently large $C = C(D)$,
\beq
\pp\left[ \sup_{t \in [t_{k}, t_{k+1}]} \left| \int_{T_1 \wedge \tau}^{t \wedge \tau} p_{s, t_k} \d N_s \right| > C \frac{ \left( \log ( \eta_{T_1} / \eta_{T_2} ) \log N \right)^{1/2}}{N \eta_{t_k}} \right] \leq N^{-D-10}.
\eeq
The claim now follows by a union bound over the intervals $[t_k, t_{k+1}]$ and \eqref{eqn:a-eta}. \qed

\section{Asymptotics of $\langle M' \rangle$} \label{a:M'-bd}

In this appendix we prove \eqref{eqn:M'-bd}.  By \cite[Eq. (A.16)]{cipolloni2023universality}, used for $z_1=z_2=z$, we have
\begin{equation}
\label{eq:mder}
\langle M_r'\rangle=-1+\frac{m^2+1-|z|^2u^2}{1+|z|^4u^4-m^4-2|z|^2u^2}=-1+\frac{1}{1-|z|^2u^2-m^2}.
\end{equation}
We now compute the leading order behavior of $1-|z|^2u^2-m^2$. First, we notice that \eqref{eq:mde} is equivalent to $-u=m^2-|z|^2u^2$. Using this relation we have
\begin{equation}
\label{eq:explicit}
1-|z|^2u^2-m^2=1+u-2|z|^2u^2=1-u+2u(1-|z|^2u)=1-u+2u(1-|z|^2+(1-u)|z|^2).
\end{equation}
Next, we notice that by \eqref{eq:mde}--\eqref{eq:mMrel} we have (recall $\rho=|\Im m|$)
\begin{equation}
\label{eq:expneed}
1-u=\frac{\eta}{\eta+\rho}=\frac{\eta}{\rho}+\mathcal{O}\left(\left(\frac{\eta}{\rho}\right)^2\right), \qquad\quad 1-|z|^2=-\frac{\eta}{\rho}+\rho^2+\mathcal{O}(\eta\rho).
\end{equation}
Plugging \eqref{eq:expneed} into \eqref{eq:explicit}, we obtain
\begin{equation}
\begin{split}
\label{eq:explden}
1-|z|^2u^2-m^2&=\frac{\eta}{\rho}+2\left(-\frac{\eta}{\rho}+\rho^2+|z|^2\frac{\eta}{\rho}\right)+\mathcal{O}\left(\left(\frac{\eta}{\rho}\right)^2+\eta\rho\right) \\
&=\frac{\eta}{\rho}+2\rho^2+\mathcal{O}\left(\left(\frac{\eta}{\rho}\right)^2+\eta\rho\right).
\end{split}
\end{equation}
Using \eqref{eq:explden} in \eqref{eq:mder}, we obtain
\[
\big|\langle M_r'\rangle\big|\le \frac{1}{2\rho_r^2+\eta_r/\rho_r}\left(1+\mathcal{O}\left(\frac{\eta_r}{\rho_r}\right)\right).
\]

\section{Linear algebra inequalities}

\bel \label{lem:schwarz-ave}
On the event that
$
|\langle G_t - M_t \rangle | \leq \rho_t,
$
we have for any matrix $B\in\cc^{2N\times 2N}$ with $\| B \| \leq 1$  and for any $p,q\ge 1$ we have
\beq \label{eqn:schwarz-ave}
|\langle G_t^p B E_i (G_t^q)^\mft E_j\rangle| \leq 2 \frac{\rho_t}{ \eta_t^{p+q-1}}.
\eeq
\eel
\proof
 By a Schwarz inequality we have
\[
\begin{split}
|\langle G_t^p B E_i (G_t^q)^\mft E_j\rangle|&\le \langle G_t^p B E_iB^*(G_t^p)^*\rangle^{1/2}\langle [(G_t^q)^\mft]^* (G_t^q)^\mft E_j\rangle^{1/2} \\
&\le \sqrt{\lVert BE_iB^*\rVert\lVert E_j\rVert} \langle G_t^p(G_t^p)^*\rangle^{1/2}\langle [(G_t^q)^\mft]^* (G_t^q)^\mft\rangle^{1/2} \\
&\le \frac{\langle \Im G_t\rangle}{\eta_t^{p+q-1}}\le 2 \frac{\rho_t}{\eta_t^{p+q-1}},
\end{split}
\]
where in the third inequality we used $\lVert E_i\rVert, \lVert E_j\rVert\le 1$.
\qed

\bel \label{lem:schwarz-iso}
On the event that
\beq
\sup_{{\bf u} , {\bf v } \in \{ {\bf x } , {\bf y } \} } | (G_t - M_t )_{ {\bf u } { \bf v}} | \leq \rho_t,
\eeq
for all choices of ${\bf u} , {\bf v } \in \{ {\bf x } , {\bf y }\}$, we have
\beq \label{eqn:schwarz-iso}
\big|(G_t^2)_{{\bf u v }}\big| \leq 2 \frac{\rho_t}{\eta_t}, \qquad\quad \big|(G_t E_i G_t^\mft E_j G_t )_{{\bf u v}}\big| \leq 2 \frac{\rho_t}{\eta_t^2}, \qquad\quad \big|(G_t G_t^\mft )_{ \bf u v }\big|\leq 2 \frac{\rho_t}{\eta_t}.
\eeq
\eel
\proof
By a Schwarz inequality we have
\[
\big|(G_t^2)_{{\bf u v }}\big|\le \sqrt{(G_tG_t^*)_{{\bf u u}}(G_t G_t^*)_{{\bf v v }}}\le \frac{1}{\eta_t}\sqrt{(\Im G_t)_{{\bf u u}}(\Im G_t)_{{\bf v v }}}\le 2 \frac{\rho_t}{\eta_t}.
\]
The estimate for $|(G_t G_t^\mft )_{ \bf u v }|$ is completely analogous and so omitted.
Finally, using again Schwarz inequality we obtain
\[
\begin{split}
\big|(G_t E_i G_t^\mft E_j G_t )_{{\bf u v}}\big|&\le \sqrt{\big[(G_t E_i G_t^\mft)(G_t E_i G_t^\mft)^*]_{{\bf u u}} (G_t^*E_jG_t)_{{\bf v v}}} \\
&\sqrt{\lVert E_iG_t^\mathfrak{t}G_tE_i\rVert\lVert e_j\rVert (G_tG_t^*)_{{\bf u u}}(G_t G_t^*)_{{\bf v v }}} \\
&\le 2 \frac{\rho_t}{\eta_t^2},
\end{split}
\]
where in the second inequality we used the norm bounds $\lVert E_i\rVert, \lVert E_j\rVert\le 1$ and $\lVert G_t\rVert\le 1/\eta_t$.

\qed

\bibliography{ev_bib}
\bibliographystyle{abbrv}

\end{document}